\documentclass[3p,times]{elsarticle}

\usepackage{lineno,hyperref}

\pdfoutput=1

\usepackage{amsmath}
\usepackage{tikz}
\usepackage{graphicx}    
\usepackage{subfig}
\usepackage{comment}
\usepackage{appendix}
\usepackage{mathtools}
\usepackage{comment}   
\usepackage{setspace}
\usepackage{algorithm}
\usepackage{amsfonts}
\usepackage[colorinlistoftodos]{todonotes}
\usepackage{algpseudocode}
\usepackage{multirow} 
\usepackage[nodots,nocompress]{numcompress}

\begin{document}
\begin{frontmatter}

\title{Alternating direction implicit time integrations for finite difference acoustic wave propagation: Parallelization and convergence}




\author[UPC]{B. Otero}
\ead{botero@ac.upc.edu}

\author[BSC,UCV]{O. Rojas}
\ead{rojasotilio@gmail.com}

\author[ETSETB]{F. Moya}
\ead{ferranmoya@gmail.com}

\author[SDSU]{J. Castillo}
\ead{jcastillo@mail.sdsu.edu}

\address[UPC]{Dpto. d'Arquitectura de Computadors, Universitat Polit\`{e}cnica de Catalunya, Barcelona, Spain}
\address[BSC]{Barcelona Supercomputing Center, Barcelona, Spain}
\address[UCV]{Escuela de Computaci\'on, Facultad de Ciencias, Universidad Central de Venezuela, Caracas, Venezuela}
\address[ETSETB]{Escola T\'{e}cnica Superior d'Enginyeria de Telecomunicaci\'o de Barcelona, Universitat Polit\`{e}cnica de Catalunya, Barcelona, Spain}
\address[SDSU]{Computational Science Research Center, San Diego State University, San Diego CA, EEUU}

\begin{abstract}
This work studies the parallelization and empirical convergence of two finite difference acoustic wave  propagation methods on 2-D rectangular grids, that use the same alternating direction implicit (ADI) time integration. This ADI integration is based on a second-order implicit Crank-Nicolson temporal discretization that is factored out by a Peaceman-Rachford decomposition of the time and space equation terms. In space, these methods highly diverge and apply different fourth-order accurate differentiation techniques. The first method uses compact finite differences (CFD) on nodal meshes that requires solving tridiagonal linear systems along each grid line, while the second one employs staggered-grid mimetic finite differences (MFD). For each method, we implement three parallel versions: (i) a multithreaded code in Octave, (ii) a C++ code that exploits OpenMP loop parallelization, and (iii) a CUDA kernel for a NVIDIA GTX 960 Maxwell card. In these implementations, the main source of parallelism is the simultaneous ADI updating of each wave field matrix, either column-wise or row-wise, according to the differentiation direction. In our numerical applications, the highest performances are displayed by the CFD and MFD CUDA codes that achieve speedups of 7.21x and 15.81x, respectively, relative to their C++ sequential counterparts with optimal compilation flags. Our test cases also allow to assess the numerical convergence and accuracy of both methods. In a problem with exact harmonic solution, both methods exhibit convergence rates close to 4 and the MDF accuracy is practically higher. Alternatively, both convergences decay to second order on smooth problems with severe gradients at boundaries, and the MDF rates degrade in highly-resolved grids leading to larger inaccuracies. This transition of empirical convergences agrees with the nominal truncation errors in space and time.   
\end{abstract}



\end{frontmatter}

\section{Introduction}

Finite difference (FD) modeling methods have been actively applied to seismic wave propagation problems on the last fifty years. In this area, advanced methods use staggered grids (SG) in combination to fourth- and higher-order FD formulas to resolve the wave equation in first order formulation, and then solving for all dependent variables~\cite{di2012new},~\cite{Etgun.Brien:2007:cma},~\cite{gra:1996:ssw},~\cite{leva:1988:ffd},~\cite{qian.shi.cui:2013:asgfd},~\cite{roj.et.al:2014:ldm}~\cite{saen.boh1:2004:fdm},~\cite{saen.gold:2000:mpe}. In particular, the velocity-pressure formulation of the acoustic wave equation allows stating free surface and rigid wall boundary conditions as simple Dirichlet conditions, and is also well suited for the implementation of absorbing boundaries. In a SG, discrete wave fields may be defined on separate nodal meshes, which are displaced by the half of grid spacing $\frac{h}{2}$ in one or more coordinate directions. On these grids, an unknown wave field is placed at the center of those it depends, and thus SG differentiation yields more accuracy in comparison to a traditional nodal FD stencil with an error dependent on full $h$.

FD SG methods mentioned above, among many other wave propagation schemes, employ explicit second-order Leapfrog time integrations that demand minimal storage for wave field histories, but strongly limit time stepping to enforce stability. Multi-step explicit time integrations with third- and higher-order accuracy have been progressively explored in search for wider Courant-Friedrichs-Lewy (CFL) stability ranges, while reducing numerical dispersion and other anomalies \cite{blanch.rober:1997:mlwc},~\cite{bohlen.witt:2016:tdvtd},~\cite{sei.symes:1995:danwp},~\cite{wicker.skama:2002:tsmemufts},\cite{zhang.zhang:2012}. On the other hand, alternating direction implicit (ADI) methods apply a semi-implicit integration to the wave equation that allows for large time stepping, and only involves previously computed wave fields. Thus, ADI methods are one-step and low memory demanding. These methods were developed by Peaceman and Rachford \cite{peace.rach:1955:nspede} to solve parabolic equations, and then used in hyperbolic problems ~\cite{ciment1975higher},~\cite{iyengar1978high},~\cite{mckee1973high}. Since then, several ADI schemes have been focused on solving the acoustic wave equation by using compact finite differences (CFD) to discretize the spatial derivatives. Recent applications are ~\cite{abdul:2015:cfdswebd},~\cite{deng2013application},~\cite{dong.dong:2016:usss},~\cite{Qinj:2009:mdi},~\cite{kim2006h},  and~\cite{liao2014dispersion}. A CFD formula linearly couples discrete evaluations of a function and its derivative of certain order at subsequent grid points, and requires solving a banded systems of linear equations (SEL) for the latter values along a whole grid line. Given this formulation, a CFD stencil presents a shorter length compared to a traditional FD with similar accuracy, and this property shorten the linear system bandwidth. For instance, a set of fourth order CFD can be applied through the solution of a tridiagonal system via a fast solver like Thomas' algorithm (TA). The literature on CFD is broad, but the seminal paper by Lele \cite{lele:1992:ffd} presents Taylor-based constructions of central and lateral formulas of several order of accuracy, and their optimization for maximum resolution and low dispersion on one-way advection problems. More recent contributions on CFD can be found on cited works above. Notice that ADI CFD methods have shown a potential for fast and accurate wave simulations, combining the large ADI CFL condition with the CFD low-dispersive propagation modeling. 

Time explicit FD wave propagation methods based on local spatial stencils are well suited for domain decomposition and code parallelization. Main reasons are their formulation on structured meshes and their low data dependency among grid sub-domains. Particularly, implementations on Graphics Processing Units (GPU) have become useful and powerful wave simulation tools in 2-D and 3-D, given that realistic applications may result computationally intensive even in two spatial dimensions ~\cite{mic.kom:2010:atd},~\cite{botero.et.al:2017:pamfd},
~\cite{rubio:2014:fds},~\cite{sudar.et.al:2016:nmsw},~\cite{zhang:2011:asg},~\cite{jun.et.al:2013:mifd}. Conventional ADI methods and CFD schemes addressing for high computational efficiency must rely on fast SEL solvers. In the case of parallel implementations, SEL solvers must be amenable to parallelism and efficiently deal with data dependencies and scalability issues. On GPUs, most of ADI methods have been developed for parabolic equations as those found in ~\cite{egloff:2011:gpuf},~\cite{stefa:2009:aadi},~\cite{wei.jang:2013:padis}, and ~\cite{zhan.jia:2013:pic}. Alternatively, CFD implementations on GPU have been mostly assembled for computational fluid dynamics based on solving the Navier-Stokes equations, as given by ~\cite{vahid.beh.et.al:2013:egpu},~\cite{mandi.mathi_2017:pms}, and~\cite{sri:2015:naec}. All these implementations with fourth-order spatial accuracy have adopted the Cyclic Reduction (CR) tridiagonal solver, and apply hardware-dependent specializations and optimizations for better performance. However, hybrid CR-TA solvers have been also proposed in \cite{giles:2014:ifds} and ~\cite{kim.wu.et.al:2011:sts}. In addition, solution methods based on approximate matrix inversion have been alternatively used by \cite{tutkun:2012:gpuaho}, and Authors in \cite{duy.chri.jac:2009:pip} and \cite{dang.chri.jac:2010:pec} developed a problem-based solution approach. Thus, last word on GPU tridiagonal solvers has not yet been said, and this motivates recent the reviews and exploratory studies in ~\cite{Hwu:2011:GCG} and ~\cite{chang:2014:gits}. 

This work builds on the contribution of Cordova et.al \cite{cor.et.al:2015:cfd}. Grid types, FD spatial discretizations, and the CN-ADI time integration used by the former methods in ~\cite{cor.et.al:2015:cfd} are the same that support our current numerical schemes. However, latter schemes differ in two main numerical aspects. First, our new CFD nodal method employs exclusively fourth order compact stencils at all grid points, as given by Lele in \cite{lele:1992:ffd} and summarized in Appendix A, whereas the former method reduces accuracy to third order at boundary closures. On wave problems exhibiting boundary layers, or any source of strong gradients near boundaries, the fully fourth order method may perform more precisely. Second, the former mimetic SG method, although globally fourth order accurate, presents a formulation in terms of compact-like operators that apply differentiation in two steps. This formulation feature has no advantage at the computation level, so we here reformulate this SG method in terms of standard mimetic operators, originally presented in \cite{cas.gro:2003:maa}-\cite{cas.hym.et:2001:fsc} and detailed in Appendix B. From the implementation point of view, this work undertakes the parallelization of both new schemes on multi-core (CPU) and many-core (GPU) architectures, and evaluate their performance to measure limiting speedup. It is worthy of mention that the literature on parallel ADI FD methods for wave propagation in acoustic or more general media is rather scarce, and we have not found schemes based on CFD or MDF, rather than the pioneer implementations by Moya F. in \cite{moya:2016:pfdm}. This deficiency highly motivates this work. 

The remainder of this paper is organized in the following way. Section 2 briefly presents the formulation of the CFD and mimetic numerical methods for acoustic wave motion, along with the CPU complexity analysis of associated pseudo-codes.  In Section 3, we describe the progressive code optimization and parallelization of our methods to improve computational performance in available CPU and GPU scenarios. Section 4 states the numerical tests used to assess accuracy and convergence of these schemes, and to compare the speedup of execution times achieved by all their parallel implementations in Section 5. Our main conclusions are finally given in Section 6.

\section{Numerical methods}
\label{sec:numer_methods}
In a medium $\Omega$ with density $\rho$ and adiabatic compression $\kappa$, the propagation of acoustic waves obyes the following system of non-homogeneous differential equations
\begin{equation}\label{eqn_acoustic}
  \left\{\begin{array}{ll}
   \frac{\partial u(\vec{x},t)}{\partial t }  = -\kappa \nabla\cdot\vec{v}(\vec{x},t) + f(\vec{x},t),\\
    \noalign{\medskip}
    \rho \frac{\partial \vec{v}(\vec{x},t)}{\partial t} = -\nabla u(\vec{x},t) , \quad \vec{x} \in \Omega, t \geq 0.  
  \end{array} \right.
\end{equation}
Above, the dependent variables are the scalar pressure $u$ and the particle velocity vector $\vec {v} = (v,w)$ , while $f$ may correspond to an additional forcing term. Appropriate initial and boundary conditions must be considered in order to find a particular solution to (\ref{eqn_acoustic}), and this velocity-pressure formulation is well suited to state boundary conditions of relevance in seismic applications. For instance, a free surface can be modeled by a homogeneous Dirichlet condition on $u$, and most absorbing conditions require the computation of $u$ along the domain boundaries. However, we present the following formulation of our numerical methods on the square domain $\Omega=[0,1]\times [0,1]$, where a Dirichlet condition $u(x,y,t) = u_0(x,y,t)$ is imposed at $(x,y)\in \partial \Omega$, and velocities at those edges must be computed during the simulation. Brief comments on the easy reformulation that also computes boundary values of $u$ are appropriately given nonetheless. Initial distributions of all three dependent variables are given at $t = 0$.   


\subsection{The nodal CFD method}
\label{sec:nodal_CFD}
To simplify notation, we use a square $N\times N$ lattice to discretize $\Omega$ for a given integer $N$, and then cell sides have a common length $h = \frac{1}{N}$. Grid values of the continuous fields $u$,$v$, and $w$ are hold by matrices $U$, $V$, $W$, respectively. The result of FD differentiation of discrete pressures are collocated on matrices $U_x$ and $U_y$, according to the direction of differentiation, and similar structures are used in the case of particle velocities. Using the CFD operators $P$ and $Q$ presented in Appendix A, pressure differentiation at all grid nodes is computed by the matrix operations
\[
   U_x P^{T} = U Q^{T}, \qquad P U_y = QU.
\]
However, pressure is known along boundary edges because of Dirichlet conditions, thus this formulation also makes use of the reduced matrix $\bar U$ that only holds interior grid values. Time updating of these inner discrete pressures is based on the discretization of momentum conservation in (\ref{eqn_acoustic}), and requires approximations to $v_x$ and $w_y$. In space, this discretization in carried out by reduced operators $\bar P$ and $\bar Q$, also introduced in Appendix A, applied to discrete velocities
\begin{equation}
	\bar{V}_x \bar{P}^T = \bar{V} \bar{Q}^T,  \qquad
	\bar{P} \bar{W}_y = \bar{Q} \bar{W}.
\label{ecuacion:VxWy}
\end{equation}  
Reduced matrices $\bar{V}$ and $\bar{W}$ are obtained from whole matrices after removing first and last rows of $V$, and first and last columns of $W$, respectively. By doing so, we avoid unnecessary velocity differentiation along boundaries.  

The application of the Alternate Direction Implicit (ADI) formulation to the acoustic model (\ref{eqn_acoustic}) starts by defining the discrete vector operators
\begin{equation}
A_{1h}\begin{bmatrix}\bar U\\ V\\ W\end{bmatrix}= -
  \begin{bmatrix} \kappa \bar{V} \bar{Q}^{T} (\bar{P}^{T})^{-1} \\ \rho^{-1} U Q^{T}(P^{T})^{-1} \\ 0\end{bmatrix}, \quad
  A_{2h}\begin{bmatrix}\bar U\\ V\\ W\end{bmatrix}=-\begin{bmatrix} \kappa \bar{P}^{-1} \bar{Q} \bar{W}\\0 \\ \rho^{-1} P^{-1}Q U\end{bmatrix},
  \label{ecuacion:A1HA2H}
\end{equation}
Next, the Crank-Nicolson discretization of time derivatives leads to the fully implicit method
\begin{align}
\left(I-\frac{\Delta t}{2}A_{1h}-\frac{\Delta t}{2}A_{2h}\right)\mathbf{U_{CN}^{m+1}}  
&=&\left(I+\frac{\Delta t}{2}A_{1h}+\frac{\Delta t}{2}A_{2h}\right)\mathbf{U_{CN}^{m}}+\frac{\Delta t}{2}\left(\mathbf{A_{F}^{m}} + \mathbf{A_{F}^{m+1}}\right),
\label{ecuacion:CN}
\end{align}
where $\Delta t$ is the time step, $\mathbf{U_{CN}}$ represents the discrete solution $[\bar U, V, W]^T$ at a given discrete time, and $\mathbf{A_{F}} = [F, 0, 0]^T$ collects the source contribution $f$ at grid interior. Coupled equations of this sort are efficiently solved after the application of an ADI splitting, that only requires the solution of banded SEL for CFD differentiation along a single coordinate direction. The Peaceman-Rachford ADI decomposition of (\ref{ecuacion:CN}), as used in \cite{cor.et.al:2015:cfd}, states a calculation in two stages
\begin{eqnarray}
\label{ecua2:PR}
\left(I-\frac{\Delta t}{2}A_{1h}\right)\mathbf{U_{CN}^{*}}&=&\left(I + \frac{\Delta t}{2}A_{2h}\right)\mathbf{U_{CN}^{m}} + \frac{\Delta t}{2} \mathbf{A_{F}}^{m}
\end{eqnarray}
\begin{eqnarray}
\label{ecua3:PR2}
\left(I - \frac{\Delta t}{2}A_{2h}\right)\mathbf{U_{CN}^{m+1}}&=&\left(I+\frac{\Delta t}{2}A_{1h}\right)\mathbf{U_{CN}^{*}} + \frac{\Delta t}{2} \mathbf{A_{F}}^{m+1},
\end{eqnarray}
each of which alternates the application of operators $A_{1h}$ and $A_{2h}$. The intermediate approximation $\mathbf{U_{CN}^{*}}$ is common on these ADI formulations, but only approximations at time levels $t = m \Delta t$ and $t = (m+1) \Delta t$ are consistent solutions to (\ref{eqn_acoustic}). According to the definition of $A_{1h}$, the first stage allows an explicit calculation of vertical velocity $W^*$ in terms of its previous values at time $t = m \Delta t$, and the differentiation of pressure $U^m$ along the $y$ direction. However, this stage also couples the calculation (row by row) of intermediate $\bar{U}^*$ and $V^*$ in the following two sets of linear systems 
\begin{equation}\label{ecua1:ADI}  
  \left\{\begin{array}{ll}
   \bar{U}^* \bar{P}^T + \frac{\Delta t}{2} \kappa \bar{V}^* \bar{Q}^{T}\,=\,\bar{U}^m \bar{P}^T - \frac{\Delta t}2 
   \left(\kappa \bar{P}^{-1} \bar{Q} W^{m}\bar{P}^T \right) + \frac{\Delta t}{2} F^{m}\bar{P}^T\\
      \noalign{\medskip}
   V^* P^T + \frac{\Delta t}{2} \rho^{-1} U^* Q^T\,=\,V^m P^T .\\
      \end{array} \right.
\end{equation}
We denote the matrices on the right hand side of above equations by $A$ and $B$, respectively, which can be computed from known solutions at $t= m\Delta t$. That is, $A = \bar{U}^m \bar{P}^T - \dfrac{\Delta t}{2}\left(\kappa \bar{P}^{-1} \bar Q \bar{W}^m - F^{m}\right){\bar P}^T$ and $B = V^m P^T$. Next, rows of intermediate pressures and horizontal velocities can be computed by forward and backward substitutions from the LU decomposition of tridiagonal $P$ and $\bar P$, on the following uncoupled approximation.
\begin{equation}\label{eq:adiSolve}
  \left\{\begin{array}{ll}
  \bar{U}_{k+1}^* {\bar P}^T= A - \dfrac{\Delta t}{2} \kappa \bar{V}_{k}^* \bar Q^T \\
      \noalign{\medskip}
        V_{k+1}^* P^T = B - \dfrac{\Delta t}{2} \rho^{-1} U_{k+1}^* Q^T  \quad , \quad V_0^* = V^m  .  \\
   \end{array} \right.
\end{equation}
\noindent
This fixed point iteration stops when $(||\bar{U}_{k+1}^* - \bar{U}_{k}^*||< \varepsilon)$ and $(|| V_{k+1}^*- V_{k}^*||< \varepsilon)$, for a given tolerance $\varepsilon > 0$. 

The second stage of the ADI decomposition (\ref{ecua3:PR2}) allows computing all wavefields at time $t = (m+1) \Delta t$ through a similar procedure. The vertical velocity $V^{m+1}$ is computed explicitly using its intermediate value $V^*$ and the approximation $U_x^*$, according to the definition of $A_{2h}$. The columns of pressure $U^{m+1}$ and velocity $W^{m+1}$ are successively computed by the fixed point ADI iteration  
\begin{align}\label{eq:adiSolve2}
\begin{dcases}
	\bar{P} U_{k+1}^{m+1} &= \quad C - \dfrac{\Delta t}{2} \kappa \bar{Q} \bar{W}_{k}^{m+1} \\		
	 P W_{k+1}^{m+1} &= \quad D - \dfrac{\Delta t}{2} \rho^{-1} Q U_{k+1}^{m+1} \\		
\end{dcases},\quad
 W_0^{m+1} = W^m.
\end{align}
where matrices $C = \bar P \bar{U}^* - \dfrac{\Delta t}{2} \bar{P} \left(\kappa \bar{V}^* \bar{Q}^T (\bar{P}^{T})^{-1} - F^{m+1}\right)$ and $D = P W^*$ are obtained from intermediate variables, as one can see. In above formulation, matrix inversion is only used for notation convenience in auxiliary matrices $A$ and $C$, but LU is actually employed to resolve the $2N$ linear systems involved in their calculations. As mentioned earlier, the additional $2N$ linear systems embedded into each fixed-point iteration in (\ref{eq:adiSolve}) and (\ref{eq:adiSolve2}), are also solved by forward/backward substitutions, using the previous LU factorization of $P$ and $\bar P$.

The above CDF formulation is amenable for including absorbing damping layers that require the computation of discrete wave fields along boundaries. By simply replacing reduced operators $\bar{P}$ and $\bar{Q}$, and reduced matrices $\bar{U}$, $\bar{V}$, and $\bar{W}$, by their full matrix versions in previous formulation steps as convenient, the boundary wave fields can be also obtained. These values can be further damped by the absorbing technique proposed in \cite{cerjan:1985}. Even more, the versatile strategy developed in \cite{sochacki:1987} for the damped wave equation can be also adopted after a straighforward reformulation of our CFD method.  

A high-level implementation of above nodal CFD method is sketched in algorithm \ref{alg:nodal_method}. This algorithm presents two sequential inner loops in lines $11$ and $15$, corresponding to the two ADI stages that finally calculate  $(\bar{U},V,W)^{m+1}$ from discrete solutions at $t = m \Delta t$. For clarity, each of these stages is broken down into a separate coding procedure given in Appendix C. At each of these procedures, there is no data dependency among the $N$ linear systems solved, either row wise or column wise. Thus, each ADI set of $N$ tridiagonal systems can be solved simultaneously by our parallel methods.   
\begin{algorithm}
\caption{The nodal CFD method}
\label{alg:nodal_method}
\label{sec:complej_nodal}
\begin{algorithmic}[1]
\small
\Require \(\mathbf{U_o}, \mathbf{V_o}, \mathbf{W_o} \in \mathbb{R}^{(N+1)\times (N+1)}\) 
\Ensure \(\mathbf{U}\in \mathbb{R}^{(N+1)\times (N+1)}, \mathbf{V}\in \mathbb{R}^{(N+1)\times (N+1)}, \mathbf{W}\in \mathbb{R}^{(N+1)\times (N+1)}\)
\Procedure{Nodal}{$\mathbf{U_o, V_o, W_o, cfl, T_{sim}}$}
	\State $ \Delta t \gets \frac{cfl}{N c_{max}}$, $I_{\Delta t}=\frac{T_{sim}}{\Delta t}$ \label{nodal_timestep}
	\State $ \mathbf{U^0} \gets \mathbf{U_o} $\Comment{Initialization}		
	\State $ \mathbf{V^0} \gets \mathbf{V_o} $
	\State $ \mathbf{W^0} \gets \mathbf{W_o} $
    \State $(\mathbf{P}, \mathbf{Q}, \mathbf{\bar{P}}, \mathbf{\bar{Q}}) \gets{Nodal\_CFD\_Operators(N)}$
	\For {$ m \in 1 \dots\ I_{\Delta t}$}
		\State $ \mathbf{A} \gets (\mathbf{U}^m(\bar{:},\bar{:}) - \frac{\Delta t}{2} (K{.*}\mathbf{\bar{P}}\setminus(\mathbf{Q}\mathbf{W}^m(:,\bar{:})) - \mathbf{F}^m) \mathbf{\bar{P}}^T$ \Comment{$A$ Computation}	\label{nodal_dense_mult}
		\State $\mathbf{B} \gets \mathbf{V}^m(\bar{:},:)\mathbf{P}^T $		
		\State $ \mathbf{W^*}(:,\bar{:}) \gets \mathbf{W}^m(:,\bar{:})-\frac{\Delta t}{2}R{.*}(\mathbf{P}\setminus\mathbf{\bar{Q}}\mathbf{U}^m(\bar{:},\bar{:}))$ 
		\State $\begin{pmatrix}\mathbf{\bar{U}} , \mathbf{\bar{V}}\end{pmatrix}^* \gets$\Call{ADI-rows}{$\mathbf{U}^m,\mathbf{V}^m,\mathbf{A},\mathbf{B}$} \label{nodal_rows}
		
		\State $\mathbf{C} \gets \mathbf{\bar{P}}(\mathbf{U}^*(\bar{:},\bar{:}) - \frac{\Delta t}{2}(K{.*}(\mathbf{V}^*(\bar{:},:)\mathbf{\bar{Q}}^T)\setminus \mathbf{\bar{P}}^T - \mathbf{F}^{m+1}) $\Comment{$C$ Computation}
		\State $ \mathbf{D} \gets \mathbf{\bar{P}}\mathbf{W}^*(:,\bar{:}) $
		\State $\mathbf{V}^{m+1}(\bar{:},:) \gets \mathbf{V}^*(\bar{:},:)-\frac{\Delta t}{2}R{.*}(\mathbf{U}^*(\bar{:},\bar{:})\mathbf{Q}^T/\mathbf{P}^T) $
				\State $\begin{pmatrix}\mathbf{\bar{U}}^{m+1} , \mathbf{\bar{W}}^{m+1}\end{pmatrix} \gets$\Call{ADI-columns}{$\mathbf{U}^*,\mathbf{W}^*,\mathbf{C},\mathbf{D}$}\label{nodal_cols}
		
		\State $ t \gets t + \Delta t$
	\EndFor
	\State \textbf{return} $ \mathbf{U},\mathbf{V},\mathbf{W}$
\EndProcedure
\normalsize
\end{algorithmic}
\end{algorithm}

\noindent

The algorithm \ref{alg:nodal_method} employs the following parameters and matrices:
\begin{itemize}
\item $U_0, V_0, W_0$: Initial conditions for pressure and velocities on the square $N\times N$ lattice that discretizes the region $[0,1]\times[0,1]$.
\item $K, R$: Matrices of grid values of material properties $\kappa$ and $\rho^{-1}$, respectively.
\item $cfl, c_{max}$: CFL stability parameter and the maximum value of wave speed $c = \sqrt{\kappa / \rho}$. 
\item $T_{sim}$: Physical time of the global simulation. 
\end{itemize}

\noindent
Moreover, we use some special index notations in this pseudo-code. Matrix \textbf{$\bar{X}$} represents the reduced form of matrix $X$ according to the formulations given in sections \ref{sec:nodal_CFD} and \ref{sec:MFD}. Matrix indexing by \textbf{$:$} refers to every element along the corresponding dimension, while \textbf{$\bar{:}$} only denotes the interior elements $2,3,\dots,N-1$. Finally, \textbf{${.*}$} represents the element wise vector/matrix multiplication, as used in Octave and MATLAB.

\subsubsection{Complexity analysis of the nodal CFD method}
As expected, the computational cost of this method is driven by the solution of the ADI embedded SEL. LU has been traditionally employed for solving tridiagonal linear systems due to its computational speed, in addition to a high numerical stability for diagonally dominant matrices (i.e., the Gaussian elimination can proceed with no pivoting). The algorithm only performs $2N$ steps to solve a given system, that results in an algorithmic complexity of $O(N)$.  In this work, we also omit any pivoting strategy, even though CFD matrices $P$ and $\bar{P}$, are not strictly diagonally dominant. Numerical results do not reflect any instability signs arising from this restraint. 

The operation complexity depends on the grid density $N$, that sets the dimension of vectors and matrices used by the nodal CFD method. In case of reduced matrices $\bar{P}$, $\bar{Q}$, $\bar{U}$, and others, $N-1 \approx N$, for $N>>1$. The total number of time iterations simply results $I_{\Delta t}=N\frac{c_{max} T_{sim}}{cfl}$, where $T_{sim}$ is the physical simulation time. In this analysis, let us consider the matrix operations in Table \ref{tab:order_matricial_operations}.   

\begin{table}[!htb]
\caption{Computational order of matrix operations carried out by both the nodal CFD and the mimetic SG methods}
\centering
\begin{tabular}[t]{lccc}
\hline
Operation & Data type & Computational & Order \\
& & cost &\\
\hline
Assignment & Dense matrix & $T_{=}$ & $N^2$\\
\hline
Addition/Substraction & Dense matrix & $T_{+}$ & $N^2$\\
\hline
 \multirow{3}{*}[2 mm]{Product} & \multirow{3}{*}[4 mm]{Scalar x Dense matrix} & \multirow{3}{*}[4 mm]{$T_{\lambda}$} & \multirow{3}{*}[4 mm]{$N^2$}\\
 \multirow{3}{*}[1 mm]{} & \multirow{3}{*}[4 mm]{Sparse x Sparse matrices} & \multirow{3}{*}[4 mm]{$T_{SS}$} & \multirow{3}{*}[4 mm]{$N$}\\
\multirow{3}{*}[1 mm]{} & \multirow{3}{*}[4 mm]{Dense x Sparse matrices} & \multirow{3}{*}[4 mm]{$T_{XS}$} & \multirow{3}{*}[4 mm]{$N^2$}\\
\hline 
SEL solutions by LU & $N$ tridiagonal systems & $T_{SEL}$ & $N^2$ \\
\hline 
Frobenius norm & Dense matrix & $T_{N}$ & $N^2$\\
\hline
\end{tabular}
\label{tab:order_matricial_operations}
\end{table} 

Table \ref{tab:CFD_operations} lists the amount of matrix operations performed by the nodal CFD method during $I_{\Delta t}$ time iterations. From there, one can easily see that the time complexity of this scheme is of polynomial order, and highly dominated by the row wise and column wise ADI cycles (see statements 11 and 15, in the algorithm \ref{alg:nodal_method}). In particular, the complexity dominant term is $T_{CFD}\approx 32*I_{\Delta t}*k_{max}*N^2 = O(N^3)$, where $k_{max}$ is the maximum number of iterations spent by each ADI stage to converge. In all numerical experiments carried out on this work, $k_{max}$ varies from 6 to 8, and thus we expect $k_{max}$ being of $O(1)$ in similar wave propagation applications.

\begin{table}[!htb]
\caption{Total amount of matrix operations carried out by the nodal CFD method}
\centering
\begin{tabular}[t]{lcccccccc}
\hline
Section & Maximum & $T_{=}$ &  $T_{+}$ & $T_{\lambda}$ & $T_{SS}$ & $T_{XS}$ & $T_{SEL}$ & $T_{N}$\\
& iterations & & & & & &\\
\hline
Initialization & 1 & 5 & & & & & & \\
\hline
$A$ computation & $I_{\Delta t}$ & 1 & & & & 2 & 1 & \\
\hline
$W^* computation$ & $I_{\Delta t}$ & 1 & 1 & 1 & & 1 & 1 & \\
\hline
ADI initialization & $I_{\Delta t}$ & 4 & 1 & 1 & & 2 &  & \\
\hline
ADI-rows & $I_{\Delta t}*k_{max}$ & 4 & 4 & 2 & & 2 & 2 & 2\\
\hline
$C$ computation & $I_{\Delta t}$ & 1 &  &  & & 2 & 1 & \\
\hline
$V^{m+1}$ computation & $I_{\Delta t}$ & 1 & 1 & 1 & & 1 & 1 & \\
\hline
ADI initialization & $I_{\Delta t}$ & 4 & 1 & 1 & & 2 &  & \\
\hline
ADI-columns & $I_{\Delta t}*k_{max}$ & 4 & 4 & 2 & & 2 & 2 & 2\\
\hline
\end{tabular}
\label{tab:CFD_operations}
\end{table} 

\subsubsection{Additional comments on the CFD implementation}
Along our optimization and parallelization process, we implement and test slightly different versions of the ADI embedded iterations (\ref{eq:adiSolve}) and (\ref{eq:adiSolve2}), in search for a higher level of parallelism. In particular, matrix  $V_{k+1}^*$ might be computed from pressure $U_{k}^*$ in the second step of (\ref{eq:adiSolve}), so this calculation can proceeds simultaneously with the first step that updates $U_{k+1}^*$. A similar strategy can be applied in (\ref{eq:adiSolve2}) for a simultaneous computation of $U_{k+1}^{m+1}$ and $W_{k+1}^{m+1}$, making this ADI formulation more suitable for parallel platforms. This simple idea is similar to the one behind the iterative Jacobi linear solver, fully parallelizable in comparison to its Gauss Seidel counterpart. However, this Jacobi-like ADI actually requires few more iterations for convergence, than the above implementation that delays $V_{k+1}^*$ computation until $U_{k+1}^*$ is available. In addition, this ADI choice also leads to a reduction of the CFL stability limit of this nodal method around $90 \%$, and this deficiency is even higher in the case of the mimetic SG method. Thus, we discard this implementation choice in current developments.   

\subsection{The staggered MFD method}
\label{sec:MFD}
In this formulation, the grid distribution of discrete wavefields takes place on a rectangular staggered grid. In analogy with the nodal mesh used by the previous CFD method, let us consider the 1-D grid of $N$ cells of length $h$ along the $x$ direction, with nodes $x_i = i*h$, and cell centers $x_{i+1/2} = \frac{x_i + x_{i+1}}{2}$ for $i = 0, \cdots, N-1$. These grid nodes are collected in vector $X_n$, while $X_{cb}$ is another vector comprising cell centers in addition to grid edges, that is $X_{cb} = (0,x_{\frac{1}{2}},...,x_{N-{\frac{1}{2}}},1)$. After defining similar grid vectors along the $y$ direction, the staggered distribution of discrete pressures correspond to the mesh locations given by $X_{cb} \otimes Y_{cb}$, where $\otimes$ represents the standard cartesian product operator used in set theory. Alternatively, grid sites of discrete horizontal and vertical velocities are given by $X_n \otimes Y_{cb}$ and $X_{cb} \otimes Y_n$, respectively. According to these staggered distributions, matrix dimensions of discrete wavefields are $U \in \mathbb{R}^{(N+2)\times (N+2)}$, $V \in \mathbb{R}^{(N+2)\times (N+1)}$ and $W \in \mathbb{R}^{(N+1) \times N}$, in contrast to the previous nodal formulation where all these matrices are $(N+1)\times (N+1)$ square. Due to the boundary conditions, we again make use of reduced matrices $\bar V$ and $\bar W$, obtained from original $V$ and $W$ after removing rows from the former and columns from the latter, as before.

The staggered FD differentiation employs the mimetic fourth-order operators $G_4$ and $D_4$ given in Appendix B, for an approximation of partial derivatives along both directions of all three wavefields
\begin{align}
	U_x = U G_4^T, &\quad U_y = G_4 U \label{eq:mimeticU}\\
	\bar{V}_x = \bar{V} D_4^T , &\quad \bar{W}_y = D_4 \bar{W}. \label{eq:mimeticV}
\end{align}
Next, we replace the nodal CFD definition of discrete vector operators $ A_{1h}$ and $A_{2h}$ by the following staggered approach

\begin{equation}
A_{1h}\begin{bmatrix}\bar U\\ V\\ W\end{bmatrix}= -
  \begin{bmatrix} \bar{V}D_4^T \\ U G_4^T \\ 0\end{bmatrix}, \quad
  A_{2h}\begin{bmatrix}\bar U\\ V\\ W\end{bmatrix}=-\begin{bmatrix} D_4 \bar{W}\\0 \\G_4 U\end{bmatrix}.
  \label{ecuacion:A1HA2HMimetic}
\end{equation}

With above definitions at hand, the formulation of this staggered MFD method proceeds by replicating the time discretization and Peaceman-Rachford two-stage ADI decomposition, applied in the previous section. A pseudo-code of this method is given by the algorithm \ref{alg:centro_method}. For completeness, the two ADI stages are broken down into separate procedures, and given by the algorithm \ref{alg:centro_adi_solver} in Appendix C. The algorithm \ref{alg:centro_method} employs the same parameters used by the nodal algorithm \ref{alg:nodal_method}, as well as similar matrix and indexing notation.

\begin{algorithm}
\caption{The staggered MFD method}
\label{alg:centro_method}
\begin{algorithmic}[1]
\small
\Require \(\mathbf{U_o}\in \mathbb{R}^{(N+2)\times (N+2)}, \mathbf{V_o}\in \mathbb{R}^{(N+2)\times N},\mathbf{W_o}\in \mathbb{R}^{N\times (N+2)}\) 
\Ensure \(\mathbf{U}\in \mathbb{R}^{(N+2)\times (N+2)}, \mathbf{V}\in \mathbb{R}^{(N+2)\times N}, \mathbf{W}\in \mathbb{R}^{N\times (N+2)}\)
\Procedure{Staggered}{$\mathbf{U_o,V_o,W_o,cfl}$}
	\State $ \Delta t \gets \frac{cfl}{N c_{max}}$, $I_{\Delta t}=\frac{T_{sim}}{\Delta t}$ \label{centro_timestep}
	\State $ \mathbf{U^0} \gets \mathbf{U_o} $\Comment{Initialization}		
	\State $ \mathbf{V^0} \gets \mathbf{V_o} $
	\State $ \mathbf{W^0} \gets \mathbf{W_o} $
	\State $ (D_4, G_4) \gets Staggered\_MFD\_Operators(N) $
	\For {$m \in 1 \dots\ I_{\Delta t}$}
		\State $ \mathbf{A} \gets \mathbf{U}^m(\bar{:},\bar{:}) - \frac{\Delta t}{2}( K{.*}\mathbf{D_4}\mathbf{W}^m(:,\bar{:}) - \mathbf{F}^m )$ \Comment{$A$ Computation}	\label{centro_dense_mult}
     	\State $\mathbf{B} \gets \mathbf{V}^m (\bar{:},:)$
	
		\State $ \mathbf{W^*}(:,\bar{:}) \gets \mathbf{W}^m(:,\bar{:})-\frac{\Delta t}{2}R{.*} (\mathbf{G_4}\mathbf{U}^m(\bar{:},\bar{:}))$
		
		\State $\begin{pmatrix}\mathbf{\bar{U}} , \mathbf{\bar{V}}\end{pmatrix}^* \gets$\Call{ADI-rows}{$\mathbf{U}^m,\mathbf{V}^m,\mathbf{A},\mathbf{B}$} \label{centro_rows}
		
		\State $\mathbf{C} \gets \mathbf{U}^*(\bar{:},\bar{:}) - \frac{\Delta t}{2}( K{.*}\mathbf{V}^*(\bar{:},:)\mathbf{D_4}^T - \mathbf{F}^{m+1} )$ \Comment{$C$ Computation}
		
	    \State $ \mathbf{D} \gets \mathbf{W}^*(:,\bar{:}) $
		
		\State $\mathbf{V}^{m+1}(\bar{:},:) \gets \mathbf{V}^*(\bar{:},:)-\frac{\Delta t}{2}R{.*} (\mathbf{U}^*(\bar{:},\bar{:})\mathbf{G_4}^T)$
		
		\State $\begin{pmatrix}\mathbf{\bar{U}}^{m+1} , \mathbf{\bar{W}}^{m+1}\end{pmatrix} \gets$\Call{ADI-columns}{$\mathbf{U}^*,\mathbf{W}^*,\mathbf{C},\mathbf{D}$}\label{centro_cols}
		
		\State $ t \gets t + \Delta t$
	\EndFor
	\State \textbf{return} $ \mathbf{U},\mathbf{V},\mathbf{W}$
\EndProcedure
\normalsize
\end{algorithmic}
\end{algorithm}

\subsubsection{Complexity analysis of the staggered MFD method}
The complexity analysis of this method considers the same matrix operations defined in Table \ref{tab:order_matricial_operations}. The new Table \ref{tab:MFD_operations} lists the total amount of these operations performed by the staggered method during $I_{\Delta t}$ time iterations. Similar to the nodal scheme, the operational cost is a $N$-dependent polynomial, with a leading term given by the solution of the ADI stages (referred in statements 11 and 15 of the algorithm \ref{alg:centro_method}), that we denote as $T_{MFD}$. Specifically, $T_{MFD}\approx 28*I_{\Delta t}*k_{max}*N^2 = O(N^3)$. Thus, the total amount of operations is less than those employed by the nodal method. However, the tridiagonal structure of embedded linear systems on the former method, and our LU based solution strategy, make equivalent the operational complexity of both methods.      

\begin{table}[!htb]
\caption{Total amount of matrix operations carried out by the staggered MFD method}
\centering
\begin{tabular}[t]{lcccccccc}
\hline
Section & Maximum & $T_{=}$ &  $T_{+}$ & $T_{\lambda}$ & $T_{SS}$ & $T_{XS}$ & $T_{SEL}$ & $T_{N}$\\
& iterations & & & & & & &\\
\hline
Initialization & 1 & 5 & & & 2 & & & \\
\hline
$A$ computation & $I_{\Delta t}$ & 1 & & & & 1 & & \\
\hline
$W^* computation$ & $I_{\Delta t}$ & 1 & 1 & 1 & & 1 & & \\
\hline
ADI initialization & $I_{\Delta t}$ & 4 & 1 & 1 & & &  & \\
\hline
ADI-rows & $I_{\Delta t}*k_{max}$ & 4 & 4 & 2 & & 2 & & 2\\
\hline
$C$ computation & $I_{\Delta t}$ & 1 & & & & 1 & & \\
\hline
$V^{m+1}$ computation & $I_{\Delta t}$ & 1 & 1 & 1 & & 1 & & \\
\hline
ADI initialization & $I_{\Delta t}$ & 1 & 1 & 1 & & &  & \\
\hline
ADI-columns & $I_{\Delta t}*k_{max}$ & 4 & 4 & 2 & & 2 & & 2\\
\hline
\end{tabular}
\label{tab:MFD_operations}
\end{table}

\section{Computational optimization}
\label{Hardware}
The implementation of the methods presented in Section \ref{sec:numer_methods} takes place on an Intel(R) Core(TM) i7-4770 CPU with $4$ cores running at $3.4$ GHz clock speed, and $8$ MB of cache memory. This architecture counts with a NVIDIA GTX 960 Maxwell GPU card and the Ubuntu Linux $16.04$ LTS operating system.

\subsection{CPU implementations}

\subsubsection{Sequential versions}
Original codes of the proposed methods were developed in Octave and used the LU decomposition of CFD matrices $P$ and $\bar{P}$ to solve ADI embedded SEL \cite{cor.et.al:2015:cfd}. To better use the capabilities of our CPU architecture and apply some code optimization techniques, we migrate these original codes to C++, and use the $4.8.4-2$ version of the \textit{gcc} compiler in this process.   



The first optimization step is based on using different compilation flags, among which the $-O3$ flag significantly improved our C++ codes, leading to an immediate performance improvement. This $-O3$ flag automatically introduces the following software optimizations: loop unrolling, function in-lining and automatic vectorization. A second set of optimizing strategies are related to matrix orientation, and its application highly affects to the new CFD code. The original code in \cite{cor.et.al:2015:cfd} solved the inherent SEL's in a vector oriented fashion, i.e.\ matrix updates and system solutions are performed a row/column at a time. To take better performance advantage of the algebra libraries, a matrix oriented approach is taken for the C++ implementation. That is, updates and system solutions proceed simultaneously on the full matrix. It can be seen that most of vector operations in both algorithms CFD and MDF decouple, so they can be performed as a single matrix operation. In particular, the CFD SEL's with multiple right hand sides and same matrix, can be solved simultaneously.

\subsubsection{Parallel versions}
For multiprocessing, we use the Octave parallel package to distribute computation along rows/columns among several processes. At first, one process is used for each row/column, but this solution performs extremely bad. Little work is assigned to each process, resulting in most of computing power being wasted in managing these processes and their communication. A second attempt consists of using a small number of processes (up to $8$ in an eight core CPU), and then lists of columns/rows are assigned to each of these processes. Performance results are slightly better, however the single-threaded implementation is still faster for all tested numerical meshes, with the grid of biggest size ($N = 1024$, as detailed in Section \ref{sec:num_tests}), as the only exception. In general, results of using this Octave parallel package are poor, because of the need of data replication as required by the multiple launched processes. Up to date, no mechanism is available for fast memory sharing/synchronization among processes in this parallel package.

Due to results discussed above, we return to the original vector orientation of this problem, as it parallelizes naturally thanks to the decoupled operations at the vector level. The partition proceeds at the row/column level, distributing each vector in a thread pool that performs the required operations. Our C++ implementations follow and exploit this data orientation. The next optimization step applied to these vector-oriented versions is loop parallelization using OpenMP. OpenMP allows us to take full advantage of Intel multi-core CPU architectures. In our OpenMP loop parallelization, we have taken into account two intrinsic aspects of our numerical discretization: (i) data dependencies, and (ii) stencil dependencies. Each parallel loop collects a common set of updating instructions applied on a particular grid zone (boundary, near boundary, or interior nodes). Scheduling of OpenMP threads is established by the flag ’auto’, which invokes as many threads as the number of available CPU cores. We observe in our experiments that this option provides the best performance for a wide range of scenarios. Increasing the number of threads over the number of cores usually entails an overhead due to communication conflicts, whereas using less threads than available cores might lead to an underuse of a multi-core platform.


\subsection{GPU implementation}
Once CPU parallelization has been appropriately exploited, we address a new GPU implementation of our numerical schemes. On the available architecture, the programming language is CUDA C based on the \textit{nvcc} compiler. In this migration process of our original sequential codes, we also use the Toolkit $7.5.17$ and functions of the cuBLAS library, and only develop specific kernels for particular functionalities not supported by cuBLAS. 


We next comment on crucial implementation details. Due to the relatively small sizes of the matrices involved in our test cases (up to $1024 \times 1024$), the vector oriented workload does not yield enough parallelism to fill the GPU. Thus, we return to the matrix approach in our GPU implementations, and instead of partitioning at the row/column level, we focus on optimizing the matrix operations themselves at the element level. This allows achieving a much higher level of parallelism, given that parallelization is performed throughout the full extent of the algorithm, and not only at the row/column solvers.\footnote{In order to preserve generality of our implementations, we do not unroll the computation chain. This could yield much higher speedups at the expense of problem specificity}.

As already mentioned, we implement some missing banded matrix functionalities in cuBLAS. For instance, the Level 3 (L3) product of a banded matrix times a dense matrix is not provided, while the Level 2 (L2) multiplication of a banded matrix times a vector is indeed available. We build the L3 operation over the L2 one by decomposing a matrix-matrix multiplication into independent matrix-vector products and accumulation operations. Our first approach, using streams to overlap the matrix-vector products, achieves a very low occupancy due to the cost of kernel launching, relative to the workload of each thread. We found that best results are obtained when the whole matrix is computed in a single launch, as the work per thread increases. Here, a higher number of threads are effectively launched, increasing the opportunities to obtain better occupancy. Same results are achieved by using stream-based distances and matrix norm computations. 

On the other hand, our GPU implementations need to apply reduction operations, like for instance, those involved in matrix norm computations. In this case, parallelization of the reduction process is performed in a hierarchical fashion. Each thread adds a pair of elements, and then half of the threads add the partial results computed. This process is repeated for $\log_2(N)$ iterations, until the whole list is reduced to a single value. Several practical considerations must be considered to achieve optimal performance: minimal divergence, array indexing to optimize data locality, unrolling of some operations, etc. Bigger reductions are computed in two steps. First, one reduces in parallel all the columns of the matrix to a vector in shared memory. Then, this vector is reduced again to a single value. Here, we have used this reduction technique to calculate the matrix Frobenius norm in our new kernels, as required to test ADI convergence according to algorithms in Appendix C.


Finally, we want to notice that the stopping criterion of ADI iterations in algorithms \ref{alg:nodal_adi_solver} and \ref{alg:centro_adi_solver} of Appendix C, imposes an important constraint for GPU performance. The number of steps required until convergence is actually limited by the constant parameter $k_{max}$, that varies from 6 to 8, in our experimental set. However, convergence on the Frobenius norm of matrices $U$, $V$, and $W$, updated by GPU processing, require that a floating point value must be transferred form device to host, for testing against $\varepsilon$ tolerance. This is wasteful as the PCIe bus is slow and has a high latency. To circumvent this problem, we study the average number of steps until convergence. In our cases, this number is roughly 6, so we thus perform this minimum number of iterations until convergence start being tested, and saving most of the unnecessary device to host transfers.

Concerning our new CUDA schemes, we want to emphasize that checking for inner ADI convergence forces the synchronization of kernel threads, slowing down the global computation, and leading to a bottleneck that hampers performance. As it was said, convergence is almost always achieved in a constant number of iterations for each test case, so we opt for executing a minimum number of iterations before any checking take place.

\section{Numerical tests}
\label{sec:num_tests}
We next present the family of numerical tests that allow assessing accuracy and convergence properties of both the CFD and MDF schemes. Same tests are used to evaluate their computation speedup gained through the optimization stages described previously, but this discussion is given in next section. Under appropriate initial and Dirichlet boundary conditions (BCs), the propagation model (\ref{eqn_acoustic}) has the following exact solution

\begin{align}
u(x,y,t)&=&\left[\Gamma\left(x^k-(1-x)^k + y^k-(1-y)^k\right) + sin\left(\frac{2\pi x}{\lambda}\right) sin\left(\frac{2\pi y}{\lambda}\right) \right] cos\left(\frac{2\pi t}{T}\right).
\label{eq:harmonic_hetero}
\end{align}

\noindent
This solution $u$ is time harmonic with period $T$ and presents a spatial behaviour controlled by the polynomial coefficient $\Gamma$, being $k$ a possitive integer, and the spatial period $\lambda$ of the sinusoidal component. These four parameters allow adjusting the difficulty of the particular problem numerically solved. In particular, the test solution can switch from a purely sinusoidal function for $\Gamma = 0$, to one with significant gradients at the boundary neighbourhood in the case of $\Gamma$ and $k$ with large magnitudes. 

We first solve the case of $\Gamma = 0$ where BCs reduce to homogeneous, consider a unitary wave speed $c = 1$ (i.e., $\kappa =\rho$), and take $\lambda = \frac{1}{4}$ and $T = \frac{1}{\sqrt{2}}$. We perform ten numerical simulations of each scheme for an increasing number of grid cells $N = 16,24, 32,48,64,96,128,256,512,1024$, so grid resolution improves with the reduction of $h = N^{-1}$. In each case, the global simulation time corresponds five time periods, namely $T_{SIM} = 5T$, and the time step $\Delta t$ is taken as the maximum permitting stable calculations, where $\Delta t = h c^{-1} cfl^{max}$. The term $cfl^{max}$ is the well known CFL bound that varies between our two FD ADI schemes, and were found in \cite{cor.et.al:2015:cfd}. Especifically, we employ the $cfl^{max}$ values of $0.91$ and $0.81$ in our CFD and MDF simulations, respectively. Next, a second set of harder tests with non homogeneous BCs are also numerically solved, where $\Gamma \neq 0$. In this case, we only inspect the representative cases of $k = 2$ and $k = 9$, and also assume $\Gamma = k$ to simplify this current parametrization.

\begin{figure}[!htbp]
\centering
\noindent
\subfloat[]{\label{fig:error_homo}\includegraphics[height=5.4cm, width=0.51\textwidth]{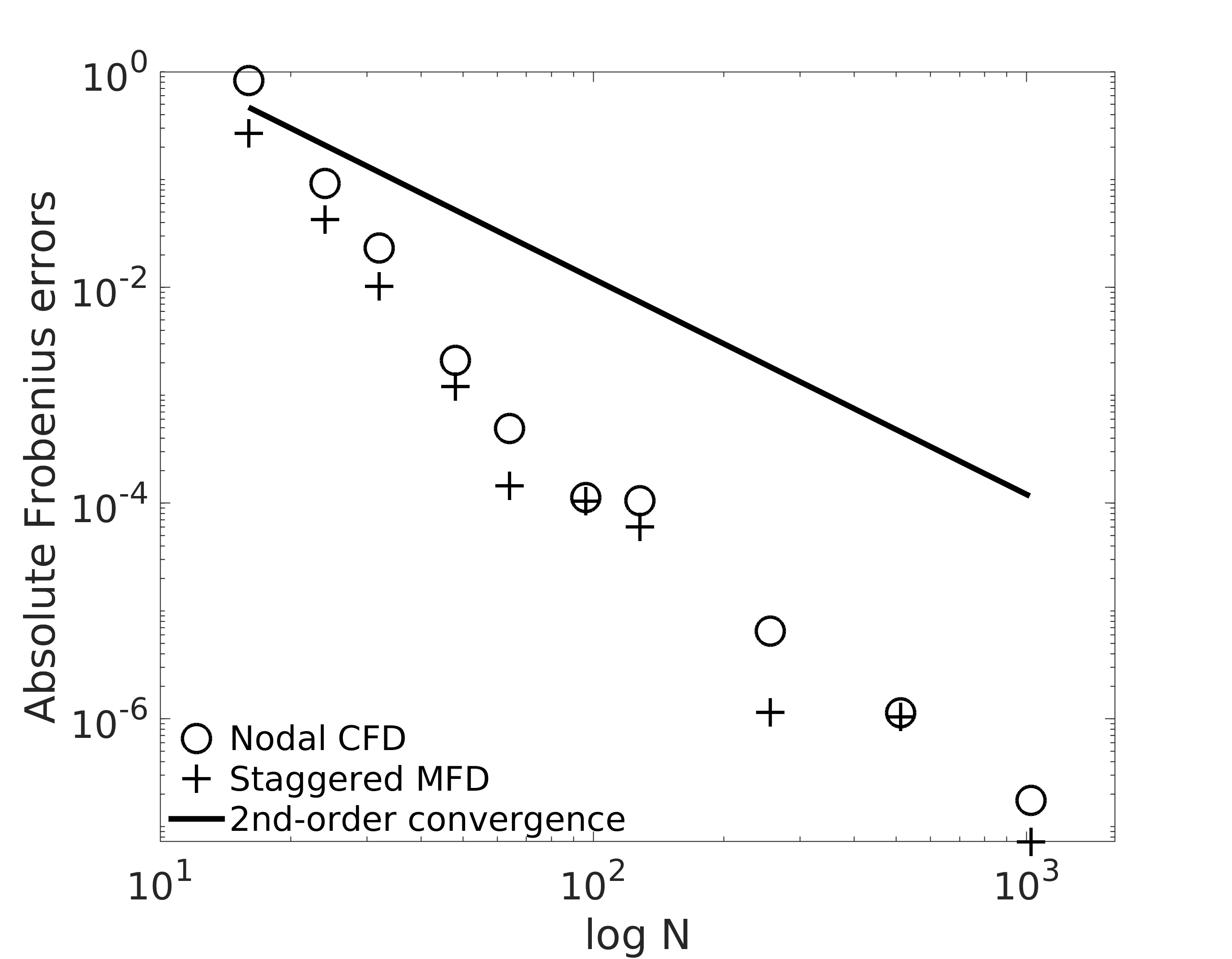}}
\subfloat[]{\label{fig:error_hetero_v2}\includegraphics[height=5.4cm, width=0.51\textwidth]{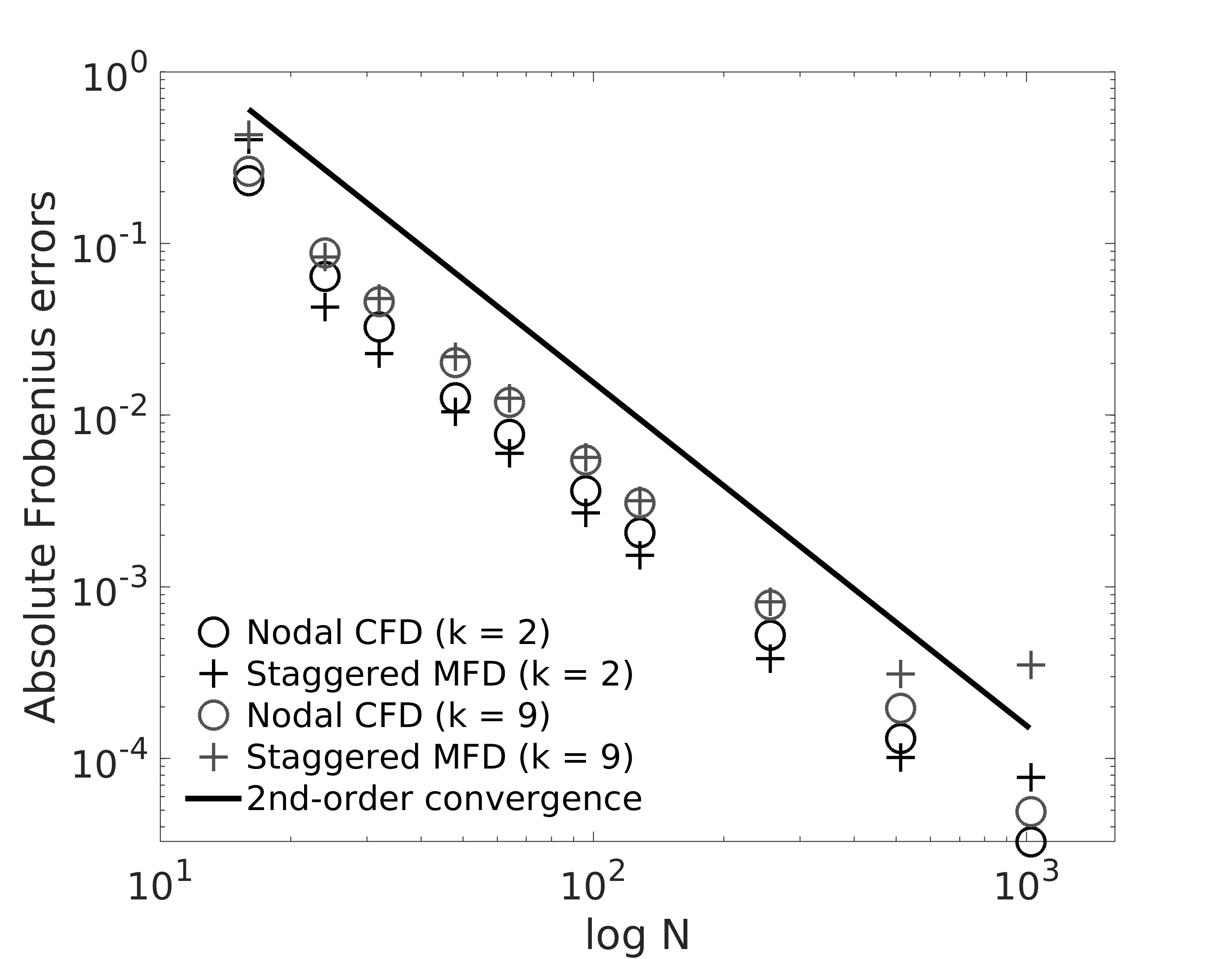}}
\caption{Comparison of absolute errors in Frobenius norm of numerical solutions to test cases under: (a) homogeneous BCs ($\Gamma = 0$); (b) heterogeneous BCs ($\Gamma = k = 2$ and $\Gamma = k = 9$)}
\label{fig:errors_num}
\end{figure}

Figure \ref{fig:errors_num} shows absolute errors in Frobenius norm found on numerical tests performed under homogeneous and heterogeneous BCs. In general, these errors steadily decay with $h$ validating our new parallel implementations. On both tests with mild difficulty, i.e., $\Gamma = 0$ and $\Gamma = k = 2$, the MFD method is slightly more accurate on most of the explored grids with $N \leq 512$, because staggering halves the grid spacing during numerical differentiation, as compared to the nodal CFD scheme. However, when solving the hardest test with $\Gamma = k = 9$, MDF errors on simulation with $N = 1024$ suffer a clear stagnation, making useless any further grid refining. Thus, we limit this experimental study, and subsequent speedup analysis, to meshes where $N \leq 1024$. 

\begin{table}[!htb]
\caption{Estimated convergence rates of the CFD simulations with time stepping given by $cfl^{MAX}$ = 0.91 }
\centering
\begin{tabular}[t]{lccc}
\hline
$N$ & $\Gamma = 0$ & $\Gamma = k = 2$ & $\Gamma = k = 9$ \\
\hline
16 &  &  & \\
\hline
24 & 5.41 & 3.17 & 2.70 \\
\hline
32 & 4.78 & 2.36 & 2.26 \\
\hline
48 & 5.92 & 2.34 & 2.02 \\
\hline
64 & 5.06 & 1.70 & 1.85 \\
\hline
96 & 3.64 & 1.87 & 1.92 \\
\hline
128 & 0.22 & 1.96 & 1.98 \\
\hline
256 & 4.02 & 1.98 & 1.98 \\
\hline
512 & 2.52 & 2.00 & 2.00 \\
\hline
1024 & 2.70 & 2.00 & 1.99\\
\hline
\bf{Average} & \bf{4.02} & \bf{2.07} & \bf{2.02} \\
\hline
\end{tabular}
\label{tab:CFD_rates}
\end{table} 

In addition, we estimate the convergence rates of the CFD and MDF solutions as $h$ reduces, and list them in tables \ref{tab:CFD_rates} and \ref{tab:MFD_rates}, respectively. These estimations are based on the standard error weighting from two consecutive simulations. For each test, we also average the resulting rates after omitting the highest and lowest values, to quantify the overall convergence behaviour avoiding main outliers. For both schemes, convergence degrades as the problem becomes harder, being almost quartic on the homogeneous test (limited by the spatial discretization), and reducing to quadratic on the hardest heterogeneous test (more constrained by the order of CN time integration). We would like to remark that accuracy and convergence results reported in \cite{cor.et.al:2015:cfd}, clearly show a MDF superiority on similar tests under homogeneous BCs. This disadvantage of the former CFD method arises from a reduction of nominal accuracy from fourth to third order at grid boundaries. Here, we opt for a fully fourth order set of CFD stencils at all grid points when implementing our current method. 


\begin{table}[!htb]
\caption{Estimated convergence rates of the MFD simulations with time stepping given by $cfl^{MAX}$ = 0.81 }
\centering
\begin{tabular}[t]{lccc}
\hline
$N$ & $\Gamma = 0$ & $\Gamma = k = 2$ & $\Gamma = k = 9$ \\
\hline
16 &  &  & \\
\hline
24 & 4.53 & 5.54 & 4.04 \\
\hline
32 & 4.94 & 2.17 & 1.93 \\
\hline
48 & 5.29 & 1.92 & 1.92 \\
\hline
64 & 7.37 & 1.94 & 1.94 \\
\hline
96 & 0.80 & 1.98 & 1.97 \\
\hline
128 & 1.94 & 1.99 & 2.01 \\
\hline
256 & 5.70 & 2.00 & 1.96 \\
\hline
512 & 0.15 & 1.91 & 1.40 \\
\hline
1024 & 3.84 & 0.40 & -0.17 \\
\hline
\bf{Average} & \bf{3.87} & \bf{1.99} & \bf{1.88} \\
\hline
\end{tabular}
\label{tab:MFD_rates}
\end{table} 


\section{Comparison of performance results}
\begin{figure}[!htbp]
\centering
\noindent
\subfloat[]{\label{fig:nodal_speed_hete}\includegraphics[height=7cm, width=0.8\textwidth]{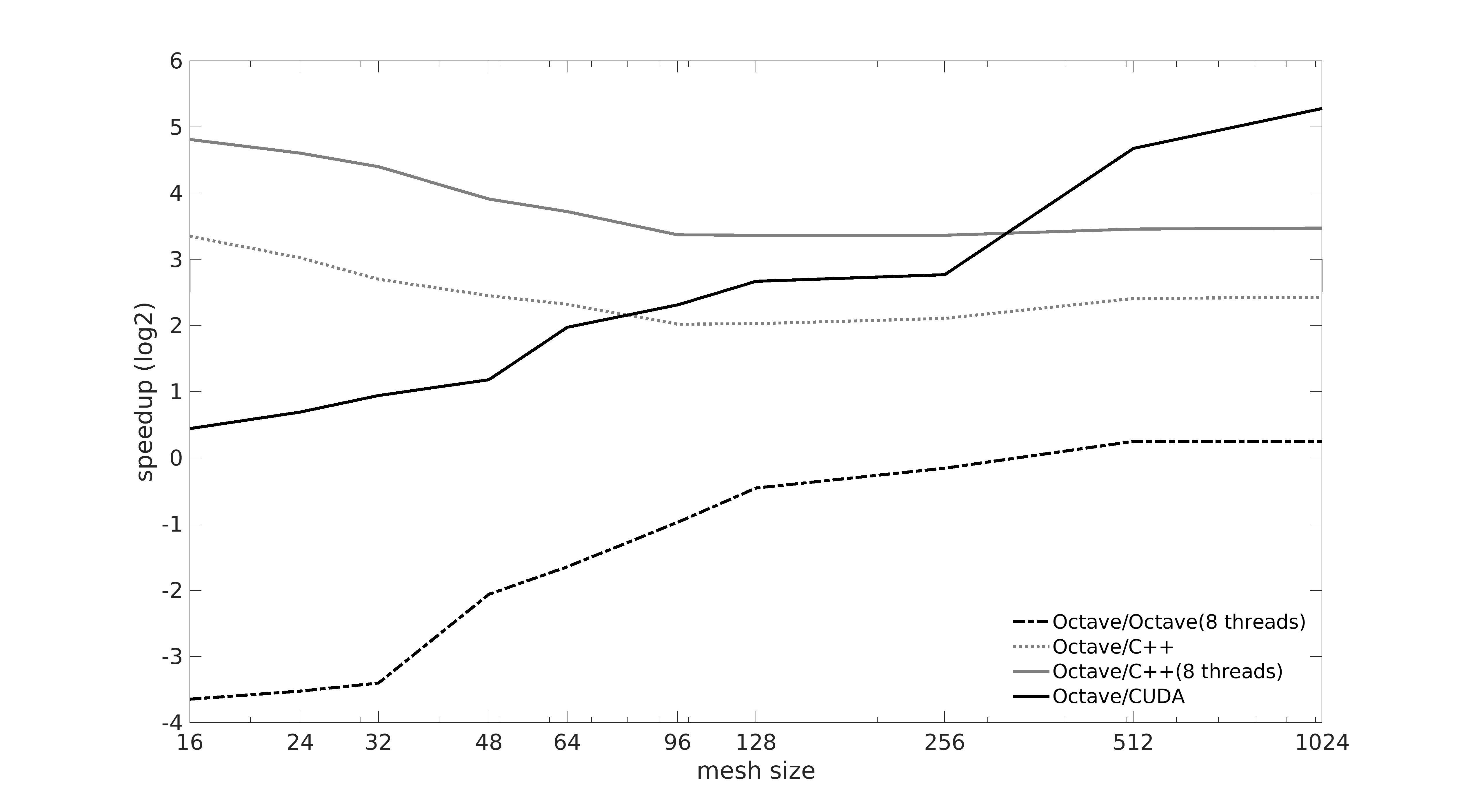}}\\

\subfloat[]{\label{fig:mimetic_speed_hete}\includegraphics[height=7cm, width=0.8\textwidth]{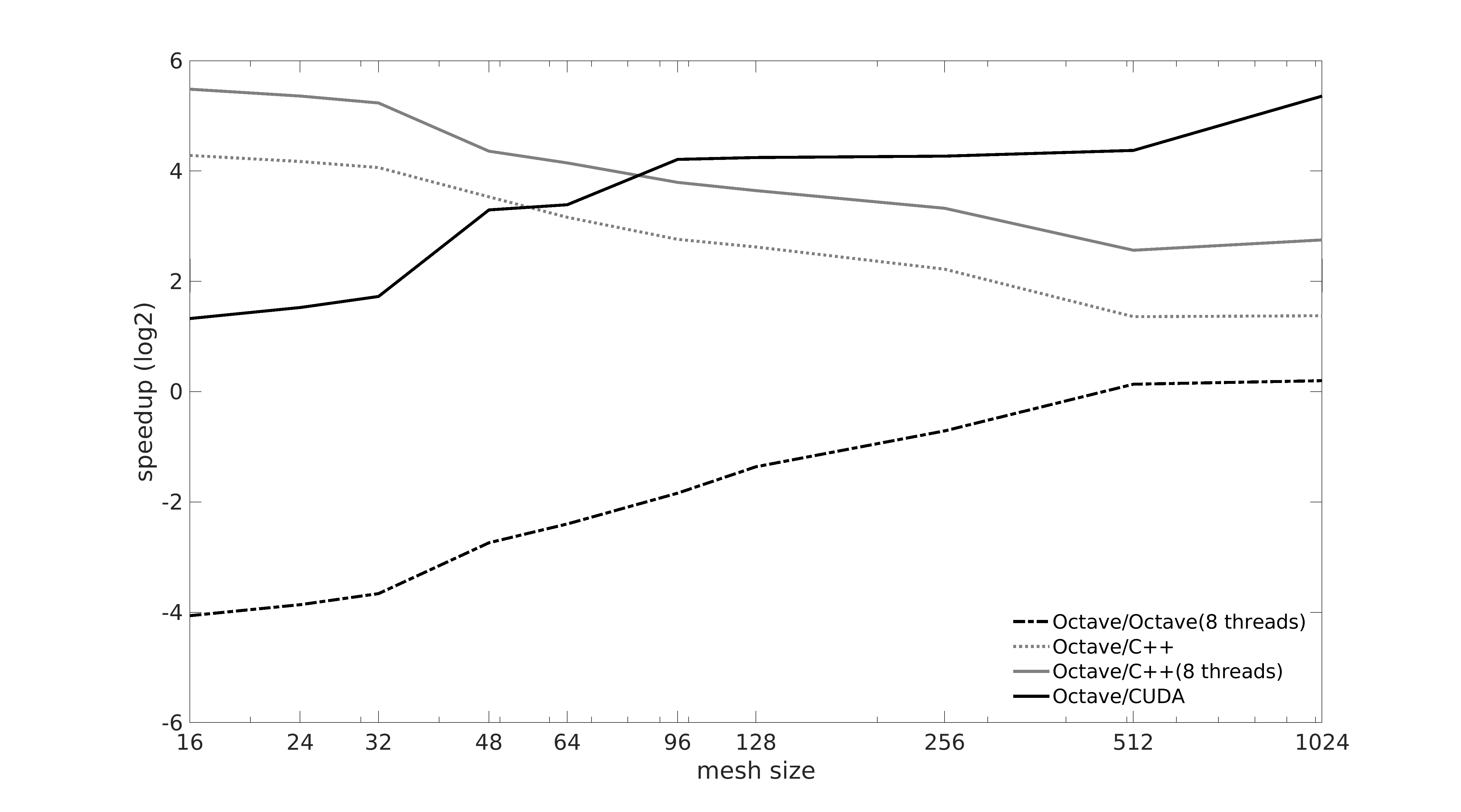}}
\caption{Speedup values in $log2$ scale achieved by the sequential C++ and parallel codes, relative to the sequential Octave implementations, of both methods the nodal CFD (a) and the staggered MFD (b), when solving the heterogeneous test for $k = 9$.}
\label{fig:errors}
\end{figure}

Figures \ref{fig:nodal_speed_hete} and \ref{fig:mimetic_speed_hete} show speedups in log2 scale, relatives to the original sequential Octave codes, of each implementation delivered by the optimization process discussed in Section \ref{Hardware}. These results were obtained after solving the numerical test under heterogeneous BCs for $k = 9$. Speedup values observed on the alternative heterogeneous or homogeneous tests, are very similar to the ones shown in these figures, confirming that evaluation of the forcing function $f$, and related operations at boundaries, minorly affect computation times.

The new CFD and MDF implementations using parallel Octave report a poor performance in most of our tests, being the former better that finally improves to a 1.19x ratio at $N$=1024. Conversely, sequential C++ implementations yield a global improvement that reaches 5.39x (CFD) and 2.6x (MFD) at the finest grid. This improvement is mainly due to software optimizations performed by the the compiler. Now, C++ parallel codes that use OpenMP achieve a higher improvement at all grids, and limiting speedups reach 11.1x (CFD) and 6.74x (MFD), at $N$ = 1024. In this case, the main reason that favors speedup is using optimized LAPACK solvers, combined to the parallelization of some operations, and to the reduction of L1-L2 cache access failures. The OpenMP implementation was almost trivial to develop from the vector oriented approach, as vectors decouple and no synchronization is needed among threads. 
Finally, we observe the best performance results on the GPU implementations, where the CFD GPU speedup grows with $N$ until 38.85x at the best resolved grid. This is a  clear indication of the high level of parallelism achieved by our GPU implementation, on the range of grid sizes tackle in this work. On the other hand, speedup of the CUDA MFD method seems to converge steadily on grids with $N \geq 96$ to a limit of $40x$, with a minor improvement at the limiting mesh, and this behavior reveals a soon high occupancy of GPUs. 

As a complementary set of speedup results, we list in Table \ref{tab:C++_speedups} only the ratios achieved by the parallel CFD and MFD methods, developed either in C++ using 8 threads or in CUDA, relative to their sequential C++ implementations. In most of our grids, these speedups improve with the mesh size. In particular, the improvements accomplished by the GPU implementations grow steadily, given the high device occupancy by our bigger matrices. It is worth reminding that grid size is limited in our tests by the numerical precision achieved by our methods, and higher refinements would not lead to more accurate results, at least for the MDF scheme. Thus, we do not test our parallel implementations with bigger matrices. Again, using the sequential C++ computations as reference, the OpenMP implementations only reach limiting improvements of 2.06x (CFD) and 2.59x (MFD) at $N=1024$, due to matrix sizes are not big enough to balance expenses on either, the parallel LAPACK LU solutions of the nodal scheme, or the several matrix-vector products of the staggered method. On the contrary, performances gained by the GPU implementations achieve maximum speedups of 7.21x (CFD) and 15.80x (MFD). Keeping matrices in GPU memory and minimizing host-device transfers were essential to obtain such performances. In these cases, the usage of streams was not effective giving that launching times were too high with respect to the kernel working times.
 
The testing process of ADI convergence (inner loops of algorithms \ref{alg:nodal_adi_solver} and \ref{alg:centro_adi_solver} in Appendix C) may result a critical restricting factor to performance, as it demands several small device-host DMA transfers. As the number of loop iterations until convergence is roughly constant, we minimize this bottleneck by delaying checking the stopping criterion until a fixed number of iterations is performed. 

\begin{table}[!htb]
\caption{Speedups achieved by the parallel CFD and MFD methods with respect to their C++ implementations in the numerical test under heterogeneous BCs ($k = 9$).}
\centering
\begin{tabular}[t]{lcccccccccc}
\hline
\textbf{Nodal CFD} & \textbf{16} & \textbf{24} & \textbf{32} & \textbf{48} & \textbf{64} & \textbf{96} & \textbf{128} & \textbf{256} & \textbf{512} & \textbf{1024} \\
\hline
C++(8 threads) & 2.75 & 2.75 & 2.76 & 2.75 & 2.64 & 2.55 & 2.52 & 2.27 & 2.07 & 2.06\\
\hline
CUDA & 0.13 & 0.18 & 0.23 & 0.41 & 0.79 & 1.22 & 1.56 & 2.54 & 4.82 & 7.21\\
\hline
\hline
\textbf{Staggered MFD} &  &  &  &  &  &  &  &  \\
\hline
C++(8 threads) & 2.29 & 2.12 & 1.96 & 1.78 & 1.98 & 2.04 & 2.03 & 2.15 & 2.31 & 2.59\\
CUDA & 0.13 & 0.29 & 0.45 & 0.85 & 1.17 & 2.73 & 3.08 & 5.94 & 8.08 & 15.81\\
\hline
\end{tabular}
\label{tab:C++_speedups}
\end{table} 

\section{Conclusions}

This paper first focuses on the computer parallelization of two fourth-order FD methods to model 2-D acoustic waves on available multi-core (CPU) and many-core (GPU) architectures. The first scheme uses CFD operators on nodal grids and requires solving a tridiagonal system of linear equations (SEL) to compute coupled derivative values along each coordinate grid line. Alternatively, the second scheme employs MFD stencils on SG for local differentiation free of SEL solution. Both methods share a Crank-Nicolson (CN) time integration combined to a Peaceman-Rachford splitting that lead to an efficient ADI updating scheme. These FD ADI methods have shown better stability and convergence properties than an explicit Leapfrog integration in a previous work \cite{cor.et.al:2015:cfd}, but the computing expenses of ADI-embedded iterations demand parallelization to be competitive on practical propagation problems. 
Along the implementation process of each method, we present a sequential C++ code with optimal compilation flags, and develop three different parallel versions. Precisely, an Octave and a C+ multithreading code, both fully exploiting array-oriented storing for matrix and vector operations, and a CUDA implementation on a NVIDIA GTX 960 Maxwell. 

Performance speedups of parallel implementations were assessed relative to the computing times of the sequential C++ code, but we additionally report results with respect to the original Octave version. These results can be insightful for someone replicating our optimization procedures for similar numerical methods. Results of GPU methods show a significant speedup gain over their parallel counterparts. Main reasons are a high occupancy of GPUs on explored grid densities, better realized by the CFD method, and a minimization of CPU-GPU transfers by delaying the convergence testing of ADI-embedded iterations, until an expected number of cycles is actually performed. Relative to Octave, the GPU CFD method achieves a speedup of 38.85x at the finest grid, while the GPU MDF reaches a higher 41x. Same performances, but measured with respect to optimal sequential C++ versions, are 7.21x for the former and 15.81x for the latter. Keeping matrices in GPU memory and minimizing host-to-device transfers result essential in obtaining such performances.       

Secondly, numerical experiments carried out in this study also enlighten convergence properties of these methods on a test suite of increasing difficulty. On a problem with an exact harmonic solution, the MDF accuracy was slightly better, and both methods exhibit convergence rates close to $4$, as the differentiation order used on the spatial discretization. However, both convergences decay to second order on smooth problems with significant gradients at boundaries, that may reveal a stronger dependence on the time integration errors. In any case, these second order experimental rates are consistent to the theoretical global convergence of $2$, and limited by the time discretization. 

\section*{Acknowledgments}
We are thankful to CAF editors and reviewers for insightful comments and revisions to this manuscript. First author was partially supported by the Generalitat de Catalunya under agreement 2017-SGR-962 and the RIS3CAT DRAC project (001-P-0 01723). The research leading to these results has received funding from the European Union’s Horizon 2020 Programme, grant agreement No. 828947, and from the Mexican Department of Energy, CONACYT-SENER Hidrocarburos grant agreement No. B-S-69926. This project has also received funding from the European Union's Horizon 2020 research and innovation programme under the Marie Sklodowska-Curie grant agreement No. 777778 (MATHROCKS). O. Rojas also thank the European Union’s Horizon 2020 Programme under the ChEESE Project, grant agreement No. 823844.

\bibliography{esqueleto}

\begin{thebibliography}{10}

\bibitem{abdul:2015:cfdswebd}
Y.~Abdulkadir.
\newblock Comparison of finite difference schemes for the wave equation based
  on dispersion.
\newblock {\em Journal of Applied Mathematics and Physics}, 3:1544--1562, 2015.

\bibitem{di2012new}
L.~Di Bartolo, C.~Dors, and W.~J. Mansur.
\newblock A new family of finite-difference schemes to solve the heterogeneous
  acoustic wave equation.
\newblock {\em Geophysics}, 77(5):T187--T199, 2012.

\bibitem{blanch.rober:1997:mlwc}
J.~O. Blanch and J.~O.~A. Robertsson.
\newblock A modified lax-wendroff correction for wave propagation in media
  described by zener elements.
\newblock {\em Geophysical Journal International}, 131(2):381--386, 1997.

\bibitem{bohlen.witt:2016:tdvtd}
T.~Bohlen and F.~Wittkamp.
\newblock Three-dimensional viscoelastic time-domain finite-difference seismic
  modelling using the staggered adams–bashforth time integrator.
\newblock {\em Geophysical Journal International}, 204(3):1781--1788, 2016.

\bibitem{cas.gro:2003:maa}
J.~E. Castillo and R.~D. Grone.
\newblock A matrix analysis approach to higher-order approximations for
  divergence and gradients satisfying a global conservation law.
\newblock {\em SIAM Journal Matrix Analysis Applications}, 25(1):128--142,
  2003.

\bibitem{cas.hym.et:2001:fsc}
J.~E. Castillo, J.~M. Hyman, M.~Shashkov, and S.~Steinberg.
\newblock Fourth- and sixth-order conservative finite difference approximations
  of the divergence and gradient.
\newblock {\em Applied Numerical Mathematics: Transactions of IMACS},
  37(1--2):171--187, 2001.

\bibitem{cerjan:1985}
C.~Cerjan, D.~Kosloff, R.~Kosloff, and M.~Reshef.
\newblock A nonreflecting boundary condition for discrete acoustic and elastic
  wave equations.
\newblock {\em Geophysics}, 50(4):705--708, 1985.

\bibitem{ciment1975higher}
M.~Ciment and S.~H. Leventhal.
\newblock Higher order compact implicit schemes for the wave equation.
\newblock {\em Mathematics of Computation}, 29(132):985--994, 1975.

\bibitem{cordova:2017:dfc}
L.~C\'ordova.
\newblock {\em Diferencias finitas compactas nodales y centro distribuidas
  aplicadas a la simulación de ondas acústicas}.
\newblock PhD thesis, Universidad Central de Venezuela, 10 2017.

\bibitem{cor.et.al:2015:cfd}
L.~C\'ordova, O.~Rojas, B.~Otero, and J.~E. Castillo.
\newblock Compact finite difference modeling of 2-d acoustic wave propagation.
\newblock {\em Journal of Computational and Applied Mathematics}, 295:83--91,
  2016.

\bibitem{duy.chri.jac:2009:pip}
D.~M. Dang, C.~C. Christara, and K.~R. Jackson.
\newblock A parallel implementation on gpus of adi finite difference methods
  for parabolic pdes with applications in finance.
\newblock {\em Canadian Applied Mathematics Quarterly}, 17(4):627--660, 2009.

\bibitem{dang.chri.jac:2010:pec}
D-M Dang, C.~C. Christara, and K.~R. Jackson.
\newblock Graphics processing unit pricing of exotic cross-currency interest
  rate derivatives with a foreign exchange volatility skew model.
\newblock {\em Concurrency and Computation Practice and Experience}, 26(9),
  2010.

\bibitem{deng2013application}
D.~Deng and C.~Zhang.
\newblock Application of a fourth-order compact adi method to solve a
  two-dimensional linear hyperbolic equation.
\newblock {\em International Journal of Computer Mathematics}, 90(2):273--291,
  2013.

\bibitem{dong.dong:2016:usss}
H.~Dongdong.
\newblock An unconditionally stable spatial sixth-order ccd-adi method for the
  two-dimensional linear telegraph equation.
\newblock {\em Numerical Algorithms}, 72(4):1103--1117, 2016.

\bibitem{egloff:2011:gpuf}
D.~Egloff.
\newblock Gpu in financial computing part iii: Adi solvers on gpus with
  application to stochastic volatility.
\newblock {\em WILMOTT Magazine}, pages 51--53, 2011.

\bibitem{vahid.beh.et.al:2013:egpu}
V.~Esfahanian, B.~Baghapour, M.~Torabzadeh, and H.~Chizari.
\newblock An efficient gpu implementation of cyclic reduction solver for
  high-order compressible viscous flow simulations.
\newblock 92, 2013.

\bibitem{Etgun.Brien:2007:cma}
J.~T. Etgen. and M.~J. O'Brien.
\newblock Computational methods for large-scale 3d acoustic finite-difference
  modeling tutorial.
\newblock {\em Geophysics}, 72(5), 2007.

\bibitem{giles:2014:ifds}
M.~Giles, E.~L\'{a}szl\'{o}, I.~Reguly, J.~Appleyard, and J.~Demouth.
\newblock Gpu implementation of finite difference solvers.
\newblock In {\em Proceedings of the 7th Workshop on High Performance
  Computational Finance}, WHPCF '14, pages 1--8, Piscataway, NJ, USA, 2014.
  IEEE Press.

\bibitem{gra:1996:ssw}
R.~W. Graves.
\newblock Simulating seismic wave propagation in 3d elastic media using
  staggered-grid finite differences.
\newblock {\em Bulletin of the Seismological Society of America},
  86(4):1091--1106, 1996.

\bibitem{Hwu:2011:GCG}
W.~M.~W. Hwu, editor.
\newblock {\em GPU Computing Gems Jade Edition}.
\newblock Morgan Kaufmann Publishers Inc., San Francisco, CA, USA, 1st edition,
  2011.

\bibitem{iyengar1978high}
S.~R. Iyengar and R.~C. Mittal.
\newblock High order difference schemes for the wave equation.
\newblock {\em International Journal for Numerical Methods in Engineering},
  12(10):1623--1628, 1978.

\bibitem{Qinj:2009:mdi}
Q.~Jinggang.
\newblock The new alternating directiion implicit difference methods for the
  wave equations.
\newblock {\em Journal of Computational and Applied Mathematics}, 230:213--223,
  2009.

\bibitem{kim.wu.et.al:2011:sts}
H.~S. Kim, S.~Wu, L.~w.~Chang, and W.~W. Hwu.
\newblock A scalable tridiagonal solver for gpus.
\newblock In {\em 2011 International Conference on Parallel Processing}, pages
  444--453, Sept 2011.

\bibitem{kim2006h}
S.~Kim and H.~Lim.
\newblock High-order schemes for acoustic waveform simulation.
\newblock {\em Applied Numerical Mathematics}, 57(4):402--414, 2007.

\bibitem{lele:1992:ffd}
S.~K. Lele.
\newblock Compact finite difference schemes with spectral-like resolution.
\newblock {\em Journal of Computational Physics}, 103:16--42, 1992.

\bibitem{leva:1988:ffd}
A.~Levander.
\newblock Fourth-order finite-difference p-sv seismograms.
\newblock {\em Geophysics}, 53:1425--1436, 1988.

\bibitem{chang:2014:gits}
C.~Li-Wen and W.~W. Hwu.
\newblock A guide for implementing tridiagonal solvers on gpus.
\newblock pages 29--44, 2014.

\bibitem{liao2014dispersion}
W.~Liao.
\newblock On the dispersion, stability and accuracy of a compact higher-order
  finite difference scheme for 3d acoustic wave equation.
\newblock {\em Journal of Computational and Applied Mathematics}, 270:571--583,
  2014.

\bibitem{mandi.mathi_2017:pms}
V.~G. Mandikas and E.~N. Mathioudakis.
\newblock A parallel multigrid solver for incompressible flows on computing
  architectures with accelerators.
\newblock {\em The Journal of Supercomputing}, 73(11):4931--4956, 2017.

\bibitem{mckee1973high}
S.~McKee.
\newblock High accuracy adi methods for hyperbolic equations with variable
  coefficients.
\newblock {\em IMA Journal of Applied Mathematics}, 11(1):105--109, 1973.

\bibitem{mic.kom:2010:atd}
D.~Michea and D.~Komatitisch.
\newblock Accelerating a three-dimensional finite-difference wave propagation
  code using gpu graphics cards.
\newblock {\em Geophysical Journal International}, 182(1):389--402, 2010.

\bibitem{moya:2016:pfdm}
F.~Moya.
\newblock Parallelization of finite difference methods: Nodal and mimetic
  solutions of the wave equation.
\newblock Master's thesis, 2016.

\bibitem{botero.et.al:2017:pamfd}
B.~Otero, J.~Franc\'es, R.~Rodr\'iguez, O.~Rojas, F.~Solano, and
  J.~Guevara-Jordan.
\newblock A performance analysis of a mimetic finite difference scheme for
  acoustic wave propagation on gpu platforms.
\newblock {\em Concurrency and Computation: Practice and Experience},
  29(4):e3880, 2017.

\bibitem{peace.rach:1955:nspede}
D.~W. Peaceman and H.~H. Rachford.
\newblock The numerical solution of parabolic and elliptic differential
  equations.
\newblock {\em Journal of the Society for Industrial and Applied Mathematics},
  3(1):28--41, 1955.

\bibitem{qian.shi.cui:2013:asgfd}
J.~Qian, S.~Wu, and R.~Cui.
\newblock Accuracy of the staggered-grid finite-difference method of the
  acoustic wave equation for marine seismic reflection modeling.
\newblock {\em Chinese Journal of Oceanology and Limnology}, 31(1):169--177,
  Jan 2013.

\bibitem{roj:2009:phd}
O.~Rojas.
\newblock {\em Modeling of rupture propagation under different friction laws
  using high-order mimetic operators}.
\newblock PhD thesis, Claremont Graduate University joint to San Diego State
  University, California, USA, 2009.

\bibitem{roj.day.cas:2008:mrp}
O.~Rojas, S.~M. Day, J.~E. Castillo, and L.~A. Dalguer.
\newblock Modelling of rupture propagation using high-order mimetic finite
  differences.
\newblock {\em Geophysical Journal International}, 172(2):631--650, 2008.

\bibitem{roj.et.al:2014:ldm}
O.~Rojas, B.~Otero, J.~E. Castillo, and S.~M. Day.
\newblock Low dispersive modeling of rayleigh waves on partly staggered grids.
\newblock {\em Computational Geosciences}, 18(1):29--43, 2014.

\bibitem{rubio:2014:fds}
F.~Rubio, M.~Hanzich, A.~Farr{\'e}s, J.~De la~Puente, and J.~M. Cela.
\newblock Finite-difference staggered grids in gpus for anisotropic elastic
  wave propagation simulation.
\newblock {\em Comput. Geosci.}, 70(C):181--189, 2014.

\bibitem{runyan2011novel}
B.~Runyan.
\newblock {\em A novel higher order finite difference time domain method based
  on the Castillo-Grone mimetic curl operator with applications concerning the
  time-dependent Maxwell equations}.
\newblock PhD thesis, San Diego State University, San Diego, USA, 2011.

\bibitem{saen.boh1:2004:fdm}
E.~H. Saenger and T.~Bohlen.
\newblock Finite-difference modeling of viscoelastic and anisotropic wave
  propagation using the rotated staggered grid.
\newblock {\em Geophysics}, 69:583--591, 2004.

\bibitem{saen.gold:2000:mpe}
E.~H. Saenger, N.~Gold, and S.~A. Shapiro.
\newblock Modeling the propagation of elastic waves using a modified
  finite-difference grid.
\newblock {\em Wave Motion}, 31:77--92, 2000.

\bibitem{sei.symes:1995:danwp}
A.~Sei and W.~Symes.
\newblock Dispersion analysis of numerical wave propagation and its
  computational consequences.
\newblock {\em Journal of Scientific Computing}, 10:1--27, 01 1995.

\bibitem{sochacki:1987}
J.~Sochacki, R.~Kubichek, J.~George, W.~R. Fletcher, and S.~Smithson.
\newblock Absorbing boundary conditions and surface waves.
\newblock {\em Geophysics}, 52(1):60--71, 1987.

\bibitem{sri:2015:naec}
A.~T. Srinath.
\newblock A novel approach to evaluating compact finite differences and similar
  tridiagonal schemes on gpu-accelerated clusters.
\newblock Master's thesis, Clemson University, 2015.

\bibitem{stefa:2009:aadi}
T.~P. Stefanski and T.~D. Drysdale.
\newblock Acceleration of the 3d adi-fdtd method using graphics processor
  units.
\newblock {\em 2009 IEEE MTT-S International Microwave Symposium Digest}, 2009.

\bibitem{sudar.et.al:2016:nmsw}
S.~Sudarmaji, S.~Sismanto, Waluyo, and B.~Soedijono.
\newblock Numerical modeling of 2d seismic wave propagation in fluid saturated
  porous media using graphics processing unit (gpu): Study case of realistic
  simple structural hydrocarbon trap.
\newblock {\em AIP Conference Proceedings}, 1755(1):100001, 2016.

\bibitem{tutkun:2012:gpuaho}
B.~Tutkun and F.~O. Edis.
\newblock A gpu application for high-order compact finite difference scheme.
\newblock {\em Computers \& Fluids}, 55:29 -- 35, 2012.

\bibitem{zhang:2011:asg}
Z.~Wang, S.~Peng, and T.~Liu.
\newblock Gpu accelerated 2-d staggered-grid finite difference seismic
  modelling.
\newblock {\em Journal of Software}, 6(8):1554--1561, 2011.

\bibitem{wei.jang:2013:padis}
Z.~Wei, B.~Jang, Y.~Zhang, and Y.~Jia.
\newblock Parallelizing alternating direction implicit solver on gpus.
\newblock {\em Procedia Computer Science}, 18:389 -- 398, 2013.

\bibitem{wicker.skama:2002:tsmemufts}
L.~J. Wicker and W.~C. Skamarock.
\newblock Time-splitting methods for elastic models using forward time schemes.
\newblock {\em Monthly Weather Review}, 130(8):2088--2097, 2002.

\bibitem{zhang.zhang:2012}
W.~Zhang, Z.~Zhang, and X.~Chen.
\newblock Three‐dimensional elastic wave numerical modelling in the presence
  of surface topography by a collocated‐grid finite ‐- difference method on
  curvilinear grids.
\newblock {\em Geophysical Journal International}, 190(1):358--378, 2012.

\bibitem{zhan.jia:2013:pic}
Y.~Zhang and Y.~Jia.
\newblock Parallelization of implicit cche2d model using cuda programming
  techniques.
\newblock In {\em World Environmental and Water Resources Congress},
  Cincinnati, Ohio, 2013.

\bibitem{jun.et.al:2013:mifd}
J.~Zhou, Y.~Cui, E.~Poyraz, D.~J. Choi, and C.~C. Guest.
\newblock Multi-gpu implementation of a 3d finite difference time domain
  earthquake code on heterogeneous supercomputers.
\newblock {\em Procedia Computer Science}, 18:1255 -- 1264, 2013.

\end{thebibliography}

\appendixpage
\begin{appendices}
\section*{Appendix A: Compact finite difference matrices}
\label{appendix_A}
Lele in \cite{lele:1992:ffd} presents a Taylor-based construction of compact finite difference (CFD) formulas of second, fourth and higher order accuracy on either nodal or staggered grids. In this work, we use the fourth order accurate CFD matrix operators $P$ and $Q$ that allow differentiation of a smooth function $v$ on $[0,1]$, whose evaluations at nodes $x_i = i*h$, for $h = \frac{1}{N}$ and $i=0,\cdots,N$, are the components of vector $V$. Nodal approximations to $\frac{dv}{dx}$ result from solving the banded linear system $P V_x = Q V$, for the following $(N+1) \times (N+1)$ matrices $P$ and $Q$  

\begin{eqnarray}
\label{matrices_P44Q44}
P_{4}&=&\left[\begin{array}{ccccccc}
\label{matriz_P44}
     6 & 18& 0&\cdots&&&0 \\
      1 & 4 & 1 & 0&\cdots&&0 \\
       0 &1&  4& 1& 0&\cdots&0 \\
     &&\ddots&\ddots&\ddots \\
      0&\cdots&&0& 1& 4&1\\
      0&\cdots&&&  0&18&6
     \end{array}\right], \\
     \label{matriz_Q44}\bigskip
Q_{4}&=&\frac{1}{h}\left[\begin{array}{ccccccc}
     -17& 9 & 9 & -1&0&\cdots \\
      -3& 0 & 3  \\
        &-3 & 0 & 3 \\
        &&\ddots&\ddots&\ddots \\
        &&   &-3& 0& 3\\
        &&  1&-9&-9& 17
     \end{array}\right].
\end{eqnarray}
\noindent
Above, matrix sub-indexes reflect the accuracy choice among the alternative CFD operators given in \cite{lele:1992:ffd}. This formulation is quite general yielding a three-parametric family of CFD operators with fourth order accuracy, at least, and those given above result from reducing matrix bandwidth to assure $P$ tridiagonally. 

An alternative pair of CFD operators can be derived from each $P$ and $Q$ set, by omitting the approximations to $\frac{dv}{dx}$ at grid edges. In this case, the reduced vector $\bar{V}_x \in {\Re}^{N-1}$ only collects values of $v$ at interior grid nodes, and the new CFD approximation to $\frac{dv}{dx}$ at point $x_1$ result from the linear combination of original formulas at nodes $x_0$ and $x_1$. The same formula is now used at $x_{N-1}$, after switching signs of CFD coefficients to account for backward differentiation. This simple process leads to the reduced matrix operators $\bar{P}$ and $\bar{Q}$ of dimensions $(N-2) \times (N-2)$ and $(N-2) \times N$, respectively, that take the following form in the case of \ref{matrices_P44Q44} given above 

\begin{eqnarray}
\label{matrices_reduced_P44Q44}
\bar{P}_{4}&=&\left[\begin{array}{ccccccc}
\label{matriz_reduced_P44}
     6 & 6& 0&\cdots&&&0 \\
      1 & 4 & 1 & 0&\cdots&&0 \\
       0 &1&  4& 1& 0&\cdots&0 \\
     &&\ddots&\ddots&\ddots \\
      0&\cdots&&0& 1& 4&1\\
      0&\cdots&&&0&6&6
     \end{array}\right], \\
     \label{matriz_reduced_Q44}\bigskip
\bar{Q}_{4}&=&\frac{1}{h}\left[\begin{array}{ccccccc}
     -1&-9 & 9 & 1&0&\cdots \\
       &-3 & 0 & 3  \\
       &   &-3 & 0 & 3 \\
       &&&\ddots&\ddots&\ddots \\
       &&& -3& 0&3 &\\
       &&& -1&-9&9 & 1
     \end{array}\right].
\end{eqnarray}

\noindent
Note that $\bar{P}_{4}$ is also tridiagonal which permits using same algorithm on the calculation of both $V$ and $\bar{V}_x$. 

\section*{Appendix B: Staggered finite differences}
\label{appendix_B}
For fourth order finite differentiation on staggered grids (SG), we employ the 1-D mimetic operators developed in \cite{cas.gro:2003:maa} and \cite{cas.hym.et:2001:fsc}, that have been applied to wave propagation by \cite{cor.et.al:2015:cfd}, \cite{roj:2009:phd},\cite{roj.day.cas:2008:mrp}, \cite{roj.et.al:2014:ldm} and \cite{runyan2011novel}. On the interval $[0,1]$, we retake the SG used in section 2.2 with $N$ cells of length $h = \frac{1}{N}$, where nodes $x_i = i*h$ and cell centers $x_{i+1/2} = \frac{x_i + x_{i+1}}{2}$ correspond to the evaluation points of two given scalar functions $v$ and $u$. Specifically, nodal $v$ evaluations are collected in vector $V \in {\Re}^{N+1}$, while evaluations of $u$ at cell centers and grid edges conform vector $U \in {\Re}^{N+2}$. Staggered differentiation is carried out by two separate FD matrix operators $D$ and $G$, where components of $DV$ approximate $\frac{dv}{dx}$ at every cell center, while $GU$ yields approximations to $\frac{du}{dx}$ at all nodes. For clarity, the following figure illustrates the spatial distributions of components of $V$, $U$, $DV$ and $GU$ on a 1-D SG
\begin{figure}[hbp]
\begin{center}
 \includegraphics[scale=0.5]{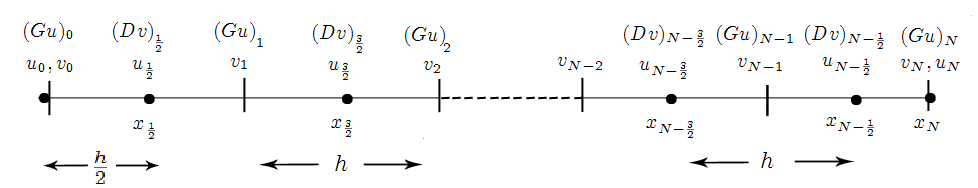}
 \label{fig:malla_encajada}
\end{center}
\end{figure}

By construction, each $D$ or $G$ operator with fourth order accuracy, comprises an individual set of SG FD stencils dependent on three free parameters $\alpha$, $\beta$ and $\gamma$. An optimal parametrization has been recently proposed in \cite{cordova:2017:dfc} to preserve the self adjointness between mimetic $D$ or $G$, an important property satisfied by their continuous differential counterparts divergence and gradient operators. However, we here adopt same parameter values $\alpha =  0$, $\beta = 0$, and $\gamma - \frac{1}{24}$ used in \cite{cor.et.al:2015:cfd}, \cite{roj.day.cas:2008:mrp}, and \cite{roj.et.al:2014:ldm} that lead to $D$ and $G$ with minimun bandwidth, to favor computational efficiency. 

\begin{eqnarray}
D_4 = \left( \begin{array}{lcccccccccl}
-\frac{4751}{5192} & \frac{909}{1298}  & \frac{6091}{15576} & -\frac{1165}{5192} &\frac{129}{2596} & -\frac{25}{15576} & 0 &\cdots \\
\,\,\,\frac{1}{24} & -\frac{9}{8}  & \frac{9}{8} & -\frac{1}{24} &0 & 0 & 0 &\cdots \\
\,\,\,\,0 & \frac{1}{24} & -\frac{9}{8}  & \frac{9}{8} & -\frac{1}{24} &0 & 0  &\cdots \\
\,\,\,\,0 & 0 &\frac{1}{24} & -\frac{9}{8}  & \frac{9}{8} & -\frac{1}{24} & 0 &\cdots \\
\end{array}
\right)
\end{eqnarray}

\begin{eqnarray}
G_4 = \left( \begin{array}{lcccccccccl}
-\frac{47888}{14245} & \frac{1790}{407}  &-\frac{14545}{9768} & \frac{8997}{16280} &-\frac{2335}{22792} &\frac{25}{9768} & 0\cdots \\
\frac{16}{105}& -\frac{31}{24} & \frac{29}{24}  & -\frac{3}{40} & \frac{1}{168} &0 & 0\cdots \\
\,\,\,0 & \frac{1}{24} & -\frac{9}{8}  & \frac{9}{8} & -\frac{1}{24} &0 & 0 \cdots \\
\,\,\,0 & 0 &\frac{1}{24} & -\frac{9}{8}  & \frac{9}{8} & -\frac{1}{24} & 0 \cdots\\
\end{array}
\right)
\end{eqnarray}

Above, matrix sub-indexes are used to reinforce the accuracy order of discretization employed in this work, because Castillo and co-workers also present operators $D$ and $G$ with alternative accuracies.
 \newpage
\section*{Appendix C: ADI solvers}
\label{appendix_C}
\begin{algorithm}
\caption{ADI iterations in the nodal CFD scheme}
\label{alg:nodal_adi_solver}
\begin{algorithmic}[1]
\small
\Require \(\mathbf{U}\in \mathbb{R}^{(N+1)\times (N+1)},\mathbf{V}\in \mathbb{R}^{(N+1)\times (N+1)},\mathbf{A}\in \mathbb{R}^{(N-2)\times (N-2)}, k >0, \varepsilon > 0 \)
\Ensure \(\mathbf{U}^*\in \mathbb{R}^{(N+1)\times (N+1)}, \mathbf{V}^*\in \mathbb{R}^{(N+1)\times (N+1)}\)
\Procedure{ADI-rows}{$\mathbf{U,V,A,B}$}
	\State $ \mathbf{U}_k \gets \mathbf{U}(\bar{:},\bar{:})$
	\State $ \mathbf{V}_k \gets \mathbf{V}(\bar{:},:)$
	\State\textbf{Repeat}\label{nodal_rsolve}
   	 \State $ \mathbf{U}_{k+1} \gets \mathbf{A} - \frac{\Delta t}{2} K{.*} \left(\mathbf{V}_k  \mathbf{\bar{Q}}^T \right)/ \mathbf{\bar{P}}^T $
   	 \State $ \mathbf{V}_{k+1} \gets \mathbf{B} - \frac{\Delta t}{2} R{.*} \left(\mathbf{U}_{k+1}  \mathbf{Q}^T \right)/ \mathbf{P}^T $
     \State $ test \gets ||\mathbf{U}_{k+1}-\mathbf{U}_{k}|| + ||\mathbf{V}_{k+1}-\mathbf{V}_{k}|| $
   	 \State $ \mathbf{U}_k \gets \mathbf{U}_{k+1}$
     \State $ \mathbf{V}_k \gets \mathbf{V}_{k+1}$
     \State $ k \gets k + 1$
  	\State \textbf{until} {$((test \leq \varepsilon)$ $||$ $(k \geq k_{max}) )$}
  	\State \textbf{return} $ \mathbf{U}^* \gets \mathbf{U}_k, \mathbf{V}^* \gets \mathbf{V}_k $
\EndProcedure

\Require \(\mathbf{U}\in \mathbb{R}^{(N+1)\times (N+1)},\mathbf{W}\in \mathbb{R}^{(N+1)\times (N+1)},\mathbf{C}\in \mathbb{R}^{(N-2)\times (N-2)}, k >0, \varepsilon > 0 \)
\Ensure \(\mathbf{U^{m+1}}\in \mathbb{R}^{(N+1)\times (N+1)}, \mathbf{W^{m+1}}\in \mathbb{R}^{(N+1)\times (N+1)}\)
\Procedure{ADI-columns}{$\mathbf{U,W,C,D}$}
	\State $ \mathbf{U}_k \gets \mathbf{U}(\bar{:},\bar{:})$
	\State $ \mathbf{W}_k \gets \mathbf{W}(:,\bar{:})$
	\State \textbf{Repeat} \label{nodal_csolve}
   	 \State $ \mathbf{U}_{k+1} \gets  \mathbf{C} - \frac{\Delta t}{2} K{.*} \mathbf{\bar{P}}  \setminus \left(\mathbf{\bar{Q}}\mathbf{W}_k \right)$
   	 \State $\mathbf{W}_{k+1} \gets  \mathbf{D} -  \frac{\Delta t}{2} R{.*} \mathbf{P}\setminus \left(\mathbf{Q}\mathbf{U}_{k+1}\right) $
     \State $ test \gets ||\mathbf{U}_{k+1}-\mathbf{U}_{k}|| + ||\mathbf{W}_{k+1}-\mathbf{W}_{k}|| $
   	 \State $ \mathbf{U}_k \gets \mathbf{U}_{k+1}$
     \State $ \mathbf{W}_k \gets \mathbf{W}_{k+1}$
     \State $ k \gets k + 1$
  	\State \textbf{until} {$((test \leq \varepsilon)$ $||$ $(k \geq k_{max}))$}
  	\State \textbf{return} $ \mathbf{U}^{m+1} \gets \mathbf{U}_k, \mathbf{W}^{m+1} \gets
  	\mathbf{W}_k $
\EndProcedure
\normalsize
\end{algorithmic}
\end{algorithm}

\begin{algorithm}
\caption{ADI iterations in the staggered MFD scheme}\label{alg:centro_adi_solver}
\begin{algorithmic}[1]
\small
\Require \(\mathbf{U}\in \mathbb{R}^{(N+2)\times (N+2)},\mathbf{V}\in \mathbb{R}^{(N+2)\times (N+1)},\mathbf{A}\in \mathbb{R}^{(N-1)\times (N-1)}, \mathbf{D_4}\in \mathbb{R}^{N\times (N+1)}, \mathbf{G_4}\in \mathbb{R}^{(N+1)\times (N+2)}, k >0, \varepsilon > 0 \)
\Ensure \(\mathbf{U^*}\in \mathbb{R}^{(N+2)\times (N+2)}, \mathbf{V^*}\in \mathbb{R}^{(N+2)\times (N+1)}\)
\Procedure{ADI-rows}{$\mathbf{U,V,A,B}$}
	\State $ \mathbf{U}_k \gets \mathbf{U}(\bar{:},\bar{:})$
	\State $ \mathbf{V}_k \gets \mathbf{V}(\bar{:},:)$
	\State \textbf{Repeat} \label{centro_rsolve}
   	 \State $ \mathbf{U}_{k+1} \gets \mathbf{A} - \frac{\Delta t}{2} K{.*} (\mathbf{V}_k{\mathbf{D}_4}^T) $ 
   	 \State $ \mathbf{V}_{k+1} \gets \mathbf{B} - \frac{\Delta t}{2} R{.*} (\mathbf{U}_{k+1}{\mathbf{G}_4}^T) $ 
     \State $ test \gets ||\mathbf{U}_{k+1}-\mathbf{U}_{k}|| + ||\mathbf{V}_{k+1}-\mathbf{V}_{k}|| $
   	 \State $ \mathbf{U}_k \gets \mathbf{U}_{k+1}$
     \State $ \mathbf{V}_k \gets \mathbf{V}_{k+1}$
     \State $ k \gets k + 1$
  	\State \textbf{until} {$((test \leq \varepsilon)$ $||$ $(k \geq k_{max}))$}
  	\State \textbf{return} $ \mathbf{U}^* \gets \mathbf{U}_k, \mathbf{V}^* \gets \mathbf{V}_k $
\EndProcedure

\Require \(\mathbf{U}\in \mathbb{R}^{(N+2)\times (N+2)},\mathbf{W}\in \mathbb{R}^{(N+1)\times (N+2)},\mathbf{C}\in \mathbb{R}^{(N-1)\times (N-1)}, k >0, \varepsilon > 0 \)
\Ensure \(\mathbf{U^{m+1}}\in \mathbb{R}^{(N+2)\times (N+2)}, \mathbf{W^{m+1}}\in \mathbb{R}^{(N+1)\times (N+2)}\)
\Procedure{ADI-columns}{$\mathbf{U,W,C,D}$}
%
	\State $ \mathbf{U}_k \gets \mathbf{U}(\bar{:},\bar{:})$
	\State $ \mathbf{W}_k \gets \mathbf{W}(:,\bar{:})$
	\State \textbf{Repeat}  \label{centro_csolve}
   	 \State $ \mathbf{U}_{k+1} \gets \mathbf{C} - \frac{\Delta t}{2} K{.*} (\mathbf{D}_4\mathbf{W}_k) $  	 
   	 \State $\mathbf{W}_{k+1} \gets \mathbf{D} - \frac{\Delta t}{2} R{.*} (\mathbf{G}_4\mathbf{U}_{k+1}) $  	
     \State $ test \gets ||\mathbf{U}_{k+1}-\mathbf{U}_{k}|| + ||\mathbf{W}_{k+1}-\mathbf{W}_{k}|| $
   	 \State $ \mathbf{U}_k \gets \mathbf{U}_{k+1}$
     \State $ \mathbf{W}_k \gets \mathbf{W}_{k+1}$
     \State $ k \gets k + 1$
  	\State \textbf{until} {$((test \leq \varepsilon)$ $||$ $(k \geq k_{max}))$}
  	\State \textbf{return} $ \mathbf{U}^{m+1} \gets \mathbf{U}_k, \mathbf{W}^{m+1} \gets
  	\mathbf{W}_k $
\EndProcedure
\normalsize
\end{algorithmic}
\end{algorithm}

\end{appendices}

\end{document}